\documentclass[reqno]{amsart}
\usepackage{amsmath, amsthm, amssymb, color}
\usepackage{graphicx}
\usepackage[mathscr]{euscript}
\usepackage{latexsym}
\usepackage{hyperref}
\hypersetup{pdfstartview={XYZ null null 1.00}} 
\usepackage{xcolor}
\usepackage{amsmath,amssymb,subeqnarray,amsfonts,graphics,epsfig,comment,color}
\usepackage{subfigure,txfonts}
\usepackage{latexsym}
 \usepackage{mathptmx}      % use Times fonts if available on your TeX system
\usepackage{color}
\usepackage[numbers,sort&compress]{natbib}
\numberwithin{equation}{section}
\newtheorem{Theorem}{Theorem}[section]
\newtheorem{Lemma}[Theorem]{Lemma}

\theoremstyle{definition}
\newtheorem{Definition}{Definition}[section]
\newtheorem{remark}[Theorem]{Remark}
\newtheorem{Proposition}[Theorem]{Proposition}

\providecommand{\norm}[1]{\left\Vert#1\right\Vert}

\newcounter{RomanNumber}

\def\be{\begin{equation}}
\def\en{\end{equation}}
\def\bs{\begin{split}}
\def\es{\end{split}}

\allowdisplaybreaks[4]
\makeatletter
\@namedef{subjclassname@2020}{\textup{2020} Mathematics Subject Classification}
\makeatother

\allowdisplaybreaks[4]

\title[Global well-posedness and instability of a generic compressible
two-fluid model]
{Global well-posedness and critical stability for a generic compressible two-fluid model}

\author{Guochun Wu}
\address{Guochun Wu \newline School of Mathematics and Statistics, Xiamen University of Technology, Xiamen, Fujian 361024, P.R. China.}
\email{guochunwu@126.com}
\author{Yunshun Wu}
\address{Yunshun Wu \newline School of Mathematical Sciences, Guizhou Normal University, Guiyang, Guizhou 550025, P.R. China.}
\email{wuyunshun1979@163.com}
\author{Yinghui Zhang}
\address{Yinghui Zhang \newline School of Mathematics and Statistics, Guangxi Normal University, Guilin, Guangxi 541004, P.R.
China} \email{yinghuizhang@mailbox.gxnu.edu.cn}

\thanks{Corresponding author: Yinghui Zhang}
\bigbreak

\date{\today}
\usepackage{hyperref}

\begin{document}
\begin{abstract}
We consider the Cauchy problem for a compressible two-fluid system with common pressure, and obtain two main results.
First, we prove global well-posedness without any low-frequency condition on the initial data: if the initial fraction densities are
close to $(1,1)$ in $H^{3}$ and the initial velocities are small in $H^{2}$, there exists a unique global strong solution. The
low-order and higher-order bounds on the gradient of the fraction densities, on the pressure and on the velocities are closed
without any time-uniform $L^{2}$ bound on the fraction densities, the low-order bound being a potential-energy equality. No decay of
the velocity or of the pressure is used, so the same argument holds in both two and three space dimensions. This removes the
low-frequency assumption in $L^{1}$ or in $\dot{B}^{s}_{1,\infty}$ required in \cite{c1, LWWZ, WYZ, WYZ2}, where the fraction
densities are recovered by integrating the decay of the first derivatives; in two dimensions that decay is only of order
$(1+t)^{-1}$ and is not integrable in time, so those methods cannot be extended to two dimensions.
Second, we identify a sharp threshold for the long-time behaviour of solutions. For $0\le s<\tfrac12$ or $\tfrac12<s\le 1$, small
initial data in $B^{s}_{1,\infty}$ give the optimal algebraic decay of the densities and of the velocities, while at the critical
index $s=\tfrac12$ the same rates hold under smallness in $\dot{B}^{1/2}_{1,1}$. This endpoint is sharp: for $\tfrac12<s\le 1$ the
fraction densities are globally algebraically unstable, although the densities and the velocities still decay algebraically.
\end{abstract}

\maketitle

%%%%%%%%%%%%%%%%%%%%%%%%%%%%%%%%%%%%%%%%%%%%%%%
\section{\leftline {\bf{Introduction.}}}
\setcounter{equation}{0}
%%%%%%%%%%%%%%%%%%%%%%%%%%%%%%%%%%%%%%%%%%%%%%%
\subsection{Background and motivation}

We consider the Cauchy problem for a compressible two-fluid model with common pressure and capillary forces. Equations of this type are used for two-phase flow in nuclear reactors, chemical processing and petroleum engineering. The derivation is given in \cite{Brennen1, Bresch1, Bresch2, Evje9, Ishii1, Prosperetti, Raja}; the modelling assumptions are set out in \cite{Bresch1}. The system reads
\begin{equation}\label{1.1}
\left\{\begin{array}{l}
\partial_{t}\bigl(\alpha^{\pm} \rho^{\pm}\bigr)+\operatorname{div}\bigl(\alpha^{\pm} \rho^{\pm} {\bf u}^{\pm}\bigr)=0, \\[1mm]
\partial_{t}\bigl(\alpha^{\pm} \rho^{\pm} {\bf u}^{\pm}\bigr)+\operatorname{div}\bigl(\alpha^{\pm} \rho^{\pm} {\bf u}^{\pm} \otimes {\bf u}^{\pm}\bigr)
+\alpha^{\pm} \nabla P
=\operatorname{div} \bigl(\alpha^{\pm} \tau^{\pm}\bigr)+\sigma^{\pm}\alpha^{\pm}
\rho^{\pm}\nabla\Delta(\alpha^{\pm} \rho^{\pm}), \\[1mm]
P=P^{\pm}\bigl(\rho^{\pm}\bigr)=A^{\pm}\bigl(\rho^{\pm}\bigr)^{\bar{\gamma}^{\pm}}.
\end{array}\right.
\end{equation}
The volume fractions satisfy \(0\le \alpha^{\pm}\le 1\) and \(\alpha^{+}+\alpha^{-}=1\). The unknowns \(\rho^{\pm}({\bf x},t)\ge 0\) and \({\bf u}^{\pm}({\bf x},t)\) are the phase densities and the velocities, and
\[
P^{\pm}\bigl(\rho^{\pm}\bigr)=A^{\pm}\bigl(\rho^{\pm}\bigr)^{\bar{\gamma}^{\pm}},\qquad \bar{\gamma}^{\pm}\ge 1,\quad A^{\pm}>0,
\]
are the pressure laws. We take \(A^{+}=A^{-}=1\). The viscous stresses are Newtonian,
\begin{equation}\label{1.2}
\tau^{\pm}:=\mu^{\pm}\bigl(\nabla {\bf u}^{\pm}+\nabla^{t}{\bf u}^{\pm}\bigr)+\lambda^{\pm}(\operatorname{div} {\bf u}^{\pm})\,\mathrm{Id},
\end{equation}
with \(\mu^{\pm}>0\) and \(2\mu^{\pm}+3\lambda^{\pm}\ge 0\), hence \(\mu^{\pm}+\lambda^{\pm}>0\). The constants \(\sigma^{\pm}\) are the capillary coefficients.
Nevertheless, the mathematical analysis of the two-fluid model presents several technical difficulties. Among these are, for instance:
\begin{itemize}
\item Viscous dissipation enters only through the terms \(\operatorname{div}(\alpha^{\pm}\tau^{\pm})\) in the momentum equations. The masses \(\alpha^{\pm}\rho^{\pm}\) are carried by the continuity equations alone, and the system is partially dissipative.

\item A phase may drop out. Once \(\alpha^{\pm}\) or \(\rho^{\pm}\) reaches zero, the corresponding mass \(\alpha^{\pm}\rho^{\pm}\) vanishes and the solution enters a single-phase region.

\item The pressure appears as the products \(\alpha^{\pm}\nabla P\). These terms lie outside conservation form, and the energy method for a single compressible fluid has to be rebuilt for the coupled phases.
\end{itemize}

On periodic domains, Bresch et al.\ \cite{Bresch1} constructed global weak solutions for \(1<\bar{\gamma}^{\pm}<6\). For the Cauchy problem in \(\mathbb{R}^{3}\), algebraic decay of strong solutions with initial data small in \(L^{1}\cap H^{3}\) was proved by Cui--Wang--Yao--Zhu \cite{c1} in the case of density-dependent viscosities, and by Li--Wang--Wu--Zhang \cite{LWWZ} for the constant viscosities \eqref{1.2}. The low-frequency threshold was then located by Wu--Yao--Zhang \cite{WYZ2}. For perturbations in \(H^{3}\cap\dot{B}^{s}_{1,\infty}\), the equilibrium is globally nonlinearly stable when \(0\le s<\frac12\), and it is unstable in the Hadamard sense when \(\frac12<s\le 1\).

In the absence of internal capillary forces (\(\sigma^{\pm}=0\)), Bresch--Huang--Li \cite{Bresch2} constructed global weak solutions in one space dimension for \(\bar{\gamma}^{\pm}>1\). In \(\mathbb{R}^{3}\), Wu--Yao--Zhang \cite{WYZ} obtained global stability for small data in \(L^{1}\cap H^{3}\), Shou--Wu--Yao--Zhang \cite{SWYZ1} treated the critical space \(\dot{B}^{-\frac12}_{2,1}\cap\dot{B}^{\frac52}_{2,1}\), and Wang--Wen \cite{WWEN} proved the stability of \cite{WYZ} under an \(L^{p}\) bound with \(1\le p<\frac65\). Wu--Wu--Zhang \cite{WWZ} obtained a critical-regularity stability theorem of the same kind as \cite{WYZ2}. For the initial-boundary value problem, Wu--Yao--Zhang \cite{WYZ2025} proved existence of a unique global strong solution for data close to equilibrium in \(H^{4}(\Omega)\). For unequal pressures satisfying
\[
P^+(\rho^+)-P^-(\rho^-)
=(\rho^+)^{\overline{\gamma}^+}-(\rho^-)^{\overline{\gamma}^-}
=f(\alpha_-\rho_-),
\]
with capillary pressure $f\in C^3([0,\infty))$, we refer to~\cite{Evje9, SWYZ2, WYZ0} and the references therein.

In summary, the partially dissipative and non-symmetric nature of the equal-pressure system \eqref{1.1} has forced all existing works \cite{c1, LWWZ, WYZ2} to impose a low-frequency smallness assumption on the initial perturbation, typically measured in $L^1$ or $\dot{B}^{s}_{1,\infty}$ with $s\in[0,1]$. This assumption is not merely technical: it yields sufficiently fast decay of the velocities, which in turn is used to control the nonlinear terms involving the non-dissipative fraction densities $\alpha^{\pm}\rho^{\pm}$. The corresponding decay rates are derived by combining energy estimates with spectral analysis. Since the fraction densities lack low-order dissipation, $\|n^{\pm}\|_{L^2}$ cannot be bounded uniformly in time (see \eqref{1.19}), and the $L^2$ estimates for the pressure and velocities in \cite{c1, LWWZ, WYZ2} can be closed only after this bound is recovered from the low-frequency decay. This raises a natural and important question: does the system \eqref{1.1} remain globally well-posed without any additional low-frequency assumption on the initial perturbation?

The main purpose of this work is to remove this low-frequency restriction. In \cite{c1, LWWZ, WYZ2}, the energy estimates are closed only after a low-frequency $L^1$ or $\dot{B}^{s}_{1,\infty}$ smallness condition produces a uniform bound on $\|n^{\pm}\|_{L^2}$. Such a bound is not available in our setting: the fraction densities have no low-order dissipation, and $\|n^{\pm}\|_{L^2}$ is not uniformly bounded in time (see \eqref{1.19}). Without the low-frequency assumption, the energy estimates can be carried out only under the \textit{a priori} assumption that, for some sufficiently small $\delta>0$,
\[
\norm{\left(\nabla n^{\pm},P-\bar{P},\mathbf{u}^{\pm}\right)}_{H^2}\le\delta.
\]

At the level of low-order energy, we introduce the potential energy densities
\[
G^{\pm}(\rho^{\pm})\triangleq\rho^{\pm}\int_{\bar{\rho}^{\pm}}^{\rho^{\pm}}\frac{P^{\pm}(s)-\bar{P}}{s^2}\,\mathrm{d}s,
\]
which lead to an energy identity independent of $\|n^{\pm}\|_{L^2}$. The viscous dissipation of $\mathbf{u}^{\pm}$ then balances the time derivative of
\[
\frac12 R^{\pm}|\mathbf{u}^{\pm}|^2+R^{\pm}G^{\pm}(\rho^{\pm})+\frac{\sigma^{\pm}}{2}|\nabla R^{\pm}|^2.
\]
Higher-order energy estimates are designed in the same spirit: they avoid $\|n^{\pm}\|_{L^2}$ altogether and close the bound for $\|(\nabla n^{\pm},P-\bar{P},\mathbf{u}^{\pm})\|_{H^2}$. The mass conservation equations, together with these uniform bounds, yield a time-growth estimate for $\|n^{\pm}\|_{L^2}$. Combined with the local existence result, a continuity argument then gives global well-posedness.

The second result is the long-time decay on the scale \(s\in[0,1]\). If the initial perturbation is small in \(B^{s}_{1,\infty}\) and \(0\le s<\frac12\) or \(\frac12<s\le 1\), the phase densities and the velocities decay at the optimal algebraic rate. At \(s=\frac12\), the same rates hold under smallness in \(\dot{B}^{1/2}_{1,1}\). This endpoint is sharp. For \(\frac12<s\le 1\), the fraction densities are globally algebraically unstable, while the phase densities and the velocities keep the algebraic decay.

\subsection{New formulation of system \eqref{1.1} and main results}
The common pressure links the two densities. We use that link to pass from \eqref{1.1} to a system in the fraction densities and the velocities, and then pose the Cauchy problem.

Differentiating $\eqref{1.1}_3$ gives
\begin{equation}\label{1.3}
{\rm d}P
  =
  s_+^2
 {\rm d}\rho^+
  =
   s_-^2
    {\rm d}\rho^-,
 \quad
  {\rm where}
   \quad
    s_\pm
     :=
      \sqrt{ \frac{{\rm d}P}{{\rm d}\rho^\pm}(\rho^\pm)}.
\end{equation}
The quantity $s_\pm$ is the sound speed of phase $\pm$.

Introduce the fraction densities, as in \cite{Bresch1},
\begin{equation}\label{1.4}
R^\pm
 =
  \alpha^\pm \rho^\pm.
\end{equation}
The constraint $\alpha^++\alpha^-=1$ converts \eqref{1.4} into
\begin{equation}\label{1.5}
{\rm d}
 \rho^+
  =
   \frac{1}{\alpha_+}
    ({\rm d}R^+
       -
        \rho^+
         {\rm d}\alpha^+),
  \quad
   {\rm d}
    \rho^-
     =
      \frac{1}{\alpha_-}
       ({\rm d}R^-
         +
          \rho^-
           {\rm d}\alpha^+).
\end{equation}
Solving \eqref{1.3} and \eqref{1.5} for ${\rm d}\alpha^+$ produces
\begin{equation}\notag
{\rm d}\alpha^+
 =
  \frac{\alpha^-s_+^2}{\alpha^-\rho^+s_+^2+\alpha^+\rho^-s_-^2}
   {\rm d}R^+
    -
     \frac{\alpha^+s_-^2}{\alpha^-\rho^+s_+^2+\alpha^+\rho^-s_-^2}
   {\rm d}R^-.
\end{equation}
Inserted into \eqref{1.5}, this identity becomes
\begin{equation}\notag
{\rm d}\rho^+
 =
  \frac{s_-^2}{\alpha^-\rho^+s_+^2+\alpha^+\rho^-s_-^2}
   (\rho^-{\rm d}R^++\rho^+{\rm d}R^-),
\end{equation}
and
\begin{equation}\notag
{\rm d}\rho^-
 =
  \frac{s_+^2}{\alpha^-\rho^+s_+^2+\alpha^+\rho^-s_-^2}
   (\rho^-{\rm d}R^++\rho^+{\rm d}R^-).
\end{equation}
Together with \eqref{1.3}, the last two relations give
\begin{equation}\label{1.6}
{\rm d}P
 =
  \mathcal{C}^2(\rho^- {\rm d} R^+ +\rho^+ {\rm d}R^-),
\end{equation}
where
\begin{equation}\notag
\mathcal{C}^2
 :=
   \frac{s_+^2s_-^2}{\alpha^-\rho^+s_+^2+\alpha^+\rho^-s_-^2},
    \quad
     {\rm and}
      \quad
       s_\pm^2
        =
        \frac{{\rm d}P(\rho^\pm)}{{\rm d}\rho^\pm}
         =
          \tilde{\gamma}^\pm \frac{P(\rho^\pm)}{\rho^\pm}.
\end{equation}

A second use of $\alpha^++\alpha^-=1$ yields
\begin{equation}\label{1.7}
\frac{R^+}{\rho^+}
 +
  \frac{R^-}{\rho^-}
   =1,
    \quad
     {\rm and \ therefore}
      \quad
       \rho^-
        =
         \frac{R^-\rho^+}{\rho^+-R^+}.
\end{equation}
In view of $\eqref{1.1}_3$, the two pressures agree if and only if
\begin{equation}\notag
\varphi(\rho^+, R^+, R^-)
 :=
  P(\rho^+)
   -
    P\left(\frac{R^-\rho^+}{\rho^+-R^+}\right)
     =
      0.
\end{equation}
Thus every pair of positive constants $\tilde R^+$, $\tilde R^-$ admits a root $\tilde\rho^+>\tilde R^+$ of $\varphi(\,\cdot\,,\tilde R^+,\tilde R^-)=0$. Differentiating in $\rho^+$ gives
\begin{equation}\notag
\frac{\partial\varphi}{\partial\rho^+}(\rho^+, R^+, R^-)
 =
  s_+^2+
   s_-^2
    \frac{R^-R^+}{(\rho^+-R^+)^2},
\end{equation}
so
\begin{equation}\notag
\frac{\partial\varphi}{\partial\rho^+}(\tilde\rho^+, \tilde R^+,
\tilde R^-)>0.
\end{equation}
The implicit function theorem and \eqref{1.4}--\eqref{1.7} therefore realise $\rho^\pm$, $\alpha^\pm$ and $\mathcal{C}$ as functions of $(R^+,R^-)$:
\begin{equation}\notag
\rho^\pm
     =
      \varrho^\pm(R^+,R^-),
\qquad
      \alpha^\pm
       =
        \alpha^\pm(R^+,R^-),
         \quad
     {\rm and \ therefore}
      \quad
      \mathcal{C}=\mathcal{C}(R^+, R^-).
\end{equation}
The algebraic details are given in \cite[pp.~614]{Bresch1}.

With these coefficients, \eqref{1.1} becomes a system for $(R^+,{\bf u}^+,R^-,{\bf u}^-)$,
\begin{equation}\label{1.8}
\left\{\begin{array}{l}
\partial_{t} R^{+}+\operatorname{div}\left(R^{+} {\bf u}^{+}\right)=0, \\
\partial_{t}\left(R^{+} {\bf u}^{+}\right)+\operatorname{div}\left(R^{+} {\bf u}^{+} \otimes {\bf u}^{+}\right)+\alpha^{+} \nabla P \\
\hspace{2.5cm}=\operatorname{div}\left\{\alpha^{+}\left[\mu^{+}\left(\nabla {\bf u}^{+}+\nabla^{t} {\bf u}^{+}\right)
+\lambda^{+} \operatorname{div} {\bf u}^{+} \operatorname{Id}\right]\right\}+\alpha^+R^+\nabla\Delta R^+, \\
\partial_{t} R^{-}+\operatorname{div}\left(R^{-} {\bf u}^{-}\right)=0, \\
\partial_{t}\left(R^{-} {\bf u}^{-}\right)+\operatorname{div}\left(R^{-} {\bf u}^{-} \otimes {\bf u}^{-}\right)+\alpha^{-} \nabla P \\
\hspace{2.5cm}=\operatorname{div}\left\{\alpha^{-}\left[\mu^{-}\left(\nabla
{\bf u}^{-}+\nabla^{t} {\bf u}^{-}\right)+\lambda^{-} \operatorname{div} {\bf u}^{-}
\operatorname{Id}\right]\right\}+\alpha^-R^-\nabla\Delta R^-.
\end{array}\right.
\end{equation}
We consider the Cauchy problem for \eqref{1.8} in $\mathbb{R}^3$, with initial data
\begin{equation}\label{1.9} (R^{+}, {\bf u}^{+},R^{-},{\bf u}^{-})({\bf x},
0)=(R_{0}^{+},  {\bf u}_{0}^{+}, R_{0}^{-}, {\bf u}_{0}^{-})({\bf x})\rightarrow({\bar{R}}^{+},  {\bf 0},{\bar{R}}^{-},  {\bf 0}) \quad
\hbox{as}\quad |{\bf x}|\rightarrow\infty,
\end{equation}
where $\bar R^{\pm}>0$ are the equilibrium fraction densities. Both are normalised to $1$.

\bigskip
\par
Before stating
our main results, let us first introduce the notations and
conventions used throughout this paper.
\(H^{k}(\mathbb{R}^{3})\) is the standard Sobolev space, with norm \(\|\cdot\|_{H^{k}}\), and \(L^{p}(\mathbb{R}^{3})\), \(1\le p\le\infty\), is the standard Lebesgue space, with norm \(\|\cdot\|_{L^{p}}\). The same symbol stands for scalar functions and for vector fields. We write \(\|(f,g)\|_{X}=\|f\|_{X}+\|g\|_{X}\). The relation \(a\lesssim b\) means \(a\le Cb\), with a constant \(C\) that depends only on the fixed parameters of the system. Derivatives are abbreviated by \(\nabla=\partial_{x}=(\partial_{1},\partial_{2},\partial_{3})\), where \(\partial_{i}=\partial_{x_{i}}\), \(\nabla_{i}=\partial_{i}\), and \(\partial_{x}^{\ell}f=\nabla^{\ell}f=\nabla(\nabla^{\ell-1}f)\).

The multiplier \(\Lambda^{s}\) is defined by
\begin{equation}
\Lambda^{s}f=\mathfrak{F}^{-1}(|{\bf\xi}|^{s}\widehat f),
\qquad s\in\mathbb{R},
\nonumber
\end{equation}
where \(\widehat f=\mathfrak{F}(f)\) is the Fourier transform of \(f\). The homogeneous Sobolev space \(\dot{H}^{s}(\mathbb{R}^{3})\) is normed by
\[
\|f\|_{\dot{H}^{s}}\overset{\triangle}{=}\|\Lambda^{s}f\|_{L^{2}}.
\]
Fix a small constant \(\eta>0\) and a radial cut-off \(\phi\in C_{0}^{\infty}(\mathbb{R}^{3}_{{\bf\xi}})\) such that \(\phi({\bf\xi})=1\) for \(|{\bf\xi}|\le\eta/2\) and \(\phi({\bf\xi})=0\) for \(|{\bf\xi}|\ge\eta\). The low- and high-frequency parts of \(f\) are
\[
f^{l}=\mathfrak{F}^{-1}[\phi({\bf\xi})\widehat f],
\qquad
f^{h}=\mathfrak{F}^{-1}[(1-\phi({\bf\xi}))\widehat f].
\]
Thus \(f=f^{l}+f^{h}\) whenever \(\widehat f\) exists.

For the dyadic decomposition, let \(\varphi\in S(\mathbb{R}^{3})\) be radial, supported in the annulus
\[
\mathcal{C}=\bigl\{\xi\in\mathbb{R}^{3}:\tfrac{3}{4}\le|\xi|\le\tfrac{8}{3}\bigr\},
\]
and normalised so that
\[
\sum_{q\in\mathbb{Z}}\varphi(2^{-q}\xi)=1
\qquad\text{for all }\xi\neq 0.
\]
The blocks are
\[
\Delta_{q}f=\mathfrak{F}^{-1}\bigl(\varphi(2^{-q}\xi)\mathfrak{F}f\bigr),
\qquad q\in\mathbb{Z}.
\]
Let
\[
\mathcal{D}(\mathbb{R}^{3})
=
\bigl\{f\in S(\mathbb{R}^{3}):D^{\alpha}\hat f(0)=0
\text{ for every multi-index }\alpha\in\mathbb{N}^{3}\bigr\},
\]
and let \(\mathcal{D}'(\mathbb{R}^{3})\) be its dual. This dual coincides with the quotient \(\mathcal{S}'(\mathbb{R}^{3})/\mathcal{P}\), where \(\mathcal{P}\) is the space of polynomials. Every \(f\in\mathcal{D}'(\mathbb{R}^{3})\) admits the homogeneous Littlewood--Paley decomposition
\[
f=\sum_{q\in\mathbb{Z}}\Delta_{q}f.
\]

\begin{Definition} Given $ s\in \mathbb{R},1\leq p,r\leq +\infty$, then the
homogeneous Besov space $\dot{B}_{p,r}^{s}(\mathbb R^3)$ is defined by
$$\dot{B}_{p,r}^{s}(\mathbb{R}^{3})=\{f\in D^{'}(\mathbb{R}^{3}):\|f\|_{\dot{B}_{p,r}^{s}}<+\infty\},$$
here
$$\norm{f}_{\dot{B}_{p,r}^{s}}\triangleq\norm{2^{qs}\|\Delta_q
f\|_{L^{p}}}_{l^{r}}.$$
\end{Definition}

\medskip
\par

\begin{Theorem}\label{1mainth}Let $\bar{\rho}^\pm=\rho^\pm(1,1)$ and $\bar P=P(1,1)$ denote equilibrium states of
$\rho^\pm$ and $P$, respectively. Suppose that $R_{0}^{\pm}-1\in H^3(\mathbb{R}^3)$ and ${\bf u}_{0}^{\pm}\in H^2(\mathbb{R}^3)$, then
there exists a constant $\delta_0$ such that if
\begin{equation}\label{1.10}
K_0:=\left\|R^\pm_0-1\right\|_{H^{3}}+\left\|{\bf u}_0^\pm
\right\|_{H^{2}} \leq \delta_0,
\end{equation}
 the Cauchy problem \eqref{1.8}--\eqref{1.9} admits a unique global strong
solution $\left(R^{+}, {\bf u}^{+}, R^{-}, {\bf u}^{-}\right)$  in $\mathbb R^3\times (0,\infty)$
satisfying,  for any $0<t<\infty$, \[
\begin{array}{l}
R^{\pm}-1 \in C^{0}\left([0, \infty) ;
H^{3}\left(\mathbb{R}^{3}\right)\right) \cap C^{1}\left([0,
\infty) ;
H^{2}\left(\mathbb{R}^{3}\right)\right), \\
{\bf u}^{\pm}\in C^{0}\left([0, \infty) ;
H^{2}\left(\mathbb{R}^{3}\right)\right) \cap C^{1}\left([0,
\infty) ; L^{2}\left(\mathbb{R}^{3}\right)\right),
\end{array}
\]
\begin{align}\label{1.11}\norm{\left(R^\pm-1\right)(t)}_{L^2}\le C(1+t)^{\frac{1}{2}}K_0,
\end{align}
and
\begin{align}\begin{split}\label{1.12}
&\norm{\left(P-\bar P,\nabla R^\pm,{\bf u}^\pm\right)(t)}_{H^{2}}+\int_0^t\left(\norm{\left(\nabla P,\nabla^2 R^\pm\right)(\tau)}_{H^2}^2+\norm{\nabla{\bf u}^\pm(\tau)}_{H^2}^2\right)\mathrm{d}\tau\leq
CK_0^2.\end{split}\end{align}
\end{Theorem}

\begin{remark}
The bound \eqref{1.11} is the essential difference from \cite{c1, LWWZ, WYZ, WYZ2}. In those works, a low-frequency \(L^1\) or \(\dot{B}^{s}_{1,\infty}\) (\(s\in[0,1]\)) smallness condition on the initial data yields a time-uniform energy bound for the fraction densities. No low-frequency assumption is imposed here, and \(\|R^{\pm}-1\|_{L^2}\) is controlled only by the growing bound \eqref{1.11}.
\end{remark}

\begin{remark}
The same argument gives global well-posedness in $\mathbb{R}^{2}$. The bounds on $(\nabla n^{\pm},P-\bar{P},{\bf u}^{\pm})$, both of low order and of higher order, close without a time-uniform $L^{2}$ bound on the fraction densities. The low-order bound is the potential-energy equality, and the higher-order bound does not use this norm. No decay of the velocity or of the pressure enters the estimate.
The results of \cite{c1, LWWZ, WYZ, WYZ2} do not extend to $\mathbb{R}^{2}$. Those arguments recover the fraction densities by integrating the decay of the first derivatives, and the decay comes from a low-frequency smallness assumption in $L^{1}$ or in $\dot{B}^{s}_{1,\infty}$ ($s\in[0,1]$). In two dimensions, even under that assumption, the optimal rate is only
\[
\|\nabla(P-\bar{P},{\bf u}^{\pm})(t)\|_{L^{2}}\lesssim(1+t)^{-1}.
\]
Since $\displaystyle\int_{1}^{\infty}(1+t)^{-1}\,\mathrm{d}t$ diverges, the decay is not integrable in time, and the energy estimates of those papers cannot be closed.
\end{remark}

The second main result is concerned with large time behavior to the Cauchy problem \eqref{1.8}--\eqref{1.9}, which reads as follows.

\smallskip
\begin{Theorem}\label{2mainth}  In addition to the conditions of Theorem \ref{1mainth}, assume further that
\begin{align}\label{1.13}\begin{split}N_0=\norm{\left(R^\pm_0-1,{\bf u}_0^\pm\right)}_{\dot B^{s}_{1,\infty}}\leq \delta_0,\ s\in[0,1]\ \text{and}\ s\neq \frac{1}{2}
\end{split}\end{align}
or
\begin{align}\label{1.14}\begin{split}N_0=\norm{\left(R^\pm_0-1,{\bf u}_0^\pm\right)}_{\dot B^{\frac{1}{2}}_{1,1}}\leq \delta_0,\ s=\frac{1}{2},
\end{split}\end{align}
then for all $t\ge 0$, the following time decay rates hold:
\begin{align}\label{1.15}\begin{split}\norm{\left(P-\bar P,\nabla R^\pm,{\bf u}^\pm\right)(t)}_{L^2}\le C(K_0,N_0)(1+t)^{-\frac{3}{4}+\frac{s}{2}},
\end{split}\end{align}
\begin{align}\label{1.16}\begin{split}\norm{\left(\nabla P,\nabla^2 R^\pm,\nabla{\bf u}^\pm\right)(t)}_{H^1}\le C(K_0,N_0)(1+t)^{-\frac{5}{4}+\frac{s}{2}},
\end{split}\end{align}
and
\begin{equation}\label{1.17}
\norm{\left(R^\pm-1\right)(t)}_{L^2}\le C(K_0,N_0)(1+t)^{-\frac{1}{4}+\frac{s}{2}}.
\end{equation}
\end{Theorem}

\begin{remark}
The constraint $s\le 1$ comes
from mathematical technique  to estimate the nonlinear terms. For example, the term  $\displaystyle\bigg(\int_0^t(1+\tau)^{-\frac{3}{4}}(1+t-\tau)^{-\frac{3}{2}+s}\mathcal E^2(t)\mathrm{d}\tau\bigg)^2$ occuring in \eqref{4.8} is bounded as follows
\begin{align}\nonumber\begin{split}\bigg(\int_0^t(1+\tau)^{-\frac{3}{4}}(1+t-\tau)^{-\frac{3}{2}+s}\mathcal E^2(t)\mathrm{d}\tau\bigg)^2\le C(1+t)^{-\frac{5}{2}+2s}\mathcal E^4(t)\le C(1+t)^{-\frac{3}{2}+s}\mathcal E^4(t).
\end{split}\end{align}
In order to guarantee the validity of the final inequality, we must impose the condition $s\le 1$.
\end{remark}

\par

The final result concerns the global algebraic instability to the Cauchy problem \eqref{1.8}--\eqref{1.9}, which is stated in the following theorem.
\par
\begin{Theorem}\label{3mainth} Let $n^\pm=R^\pm-1$. In addition to the conditions of Theorem \ref{1mainth} and Theorem \ref{2mainth}, assume further that
$\frac{1}{2}<s\le 1$,  assume that the Fourier transform $\left(\widehat{ n^+_0},\widehat{ {\bf u}^+_0},\widehat{ n^-_0},\widehat{ {\bf u}^-_0}\right)$ satisfies
\begin{equation}\label{1.18}\quad\quad\widehat{ n^+_0}(\xi)=
\wedge^{-1}\operatorname{div}\
\widehat{ {\bf u}^+_0}(\xi)=\widehat{ n^-_0}(\xi)=0,~\text{\rm and }~
\wedge^{-1}\text{\rm div}\ \widehat{ {\bf u}^-_0}(\xi)\ge c_0
|\xi|^s,
\end{equation}
where  $c_0$ is a positive constant. Then  the steady state
$(1,  {\bf 0},1,  {\bf 0})$  of the Cauchy problem \eqref{1.8}--\eqref{1.9} is global algebraic unstable, that is, there exists a positive constant $c$ independent of time such
that for any large enough $t$, the global strong solution $\left(R^{+}, {\bf u}^{+}, R^{-}, {\bf u}^{-}\right)$ obtained in Theorem \ref{1mainth} satisfying
\begin{align}\label{1.19}\min\left\{\norm{{n}^{+}(t)}_{L^2},
\norm{{{n}}^{-}(t)}_{L^2}\right\}\ge c(1+t)^{-\frac{1}{4}+\frac{s}{2}}.
\end{align}
\end{Theorem}

\par
\par
\indent Now, let us illustrate the main difficulties encountered in proving Theorem \ref{1mainth}, Theorem \ref{2mainth} and Theorem \ref{3mainth}, and explain our strategies to overcome them.

For the proof of Theorem \ref{1mainth}, no low-frequency \(L^1\) or \(\dot{B}^{s}_{1,\infty}\) (\(s\in[0,1]\)) smallness is imposed on the initial data. In \cite{c1, LWWZ, WYZ, WYZ2}, that smallness produces a time-uniform bound on the fraction densities, and the estimates of the pressure and the velocities are closed only after this bound is known. The same bound is not available here. Indeed, \eqref{1.12} controls \((P-\bar P,\nabla R^{\pm},{\bf u}^{\pm})\) in \(H^2\) and the time integrals of \(\nabla P\), \(\nabla^2 n^{\pm}\) and \(\nabla{\bf u}^{\pm}\), but it does not control \(\|R^{\pm}-1\|_{L^2}\). The energy estimates are therefore made under the {\rm a priori} assumption
\[
\|(\nabla n^{\pm},P-\bar P,{\bf u}^{\pm})\|_{H^2}\le\delta,
\]
with \(\delta>0\) sufficiently small.
For the low-order energy, we introduce the potential energy densities
\[
G^{\pm}(\rho^{\pm})\triangleq\rho^{\pm}\int_{\bar{\rho}^{\pm}}^{\rho^{\pm}}\frac{P^{\pm}(s)-\bar P}{s^2}\,\mathrm{d}s.
\]
The resulting energy equality is independent of \(\|n^{\pm}\|_{L^2}\): \(R^{\pm}G^{\pm}(\rho^{\pm})\) bounds the pressure, the capillary energy bounds \(\nabla n^{\pm}\), and the dissipation comes from \({\bf u}^{\pm}\). In the higher-order estimates, \(\|n^{\pm}\|_{L^2}\) is avoided as well, and the uniform bound \eqref{1.12} is closed. The mass conservation equations and \eqref{1.12} then give the time-growth bound \eqref{1.11}. By the Cauchy--Schwarz inequality,
\[
\|(R^{\pm}-1)(t)\|_{L^2}\le C(1+t)^{\frac12}K_0.
\]
This norm may grow, and it is kept out of the {\rm a priori} assumption. Local existence and a continuity argument extend the solution to all \(t\ge 0\), since the quantity controlling the continuation is the left-hand side of \eqref{1.12}. The growth \eqref{1.11} is the essential difference from \cite{c1, LWWZ, WYZ, WYZ2}, where the corresponding norm remains bounded.

For the proof of Theorem \ref{2mainth}, the decay rates use low-frequency information that is absent from Theorem \ref{1mainth}.  The smallness \(N_0\) is assumed  for \((R^\pm_0-1,{\bf u}_0^{\pm})\): in \(\dot{B}^{s}_{1,\infty}\) when \(s\in[0,1]\) and \(s\neq\frac12\), as in \eqref{1.13}, and in \(\dot{B}^{1/2}_{1,1}\) when \(s=\frac12\), as in \eqref{1.14}. The index \(s=\frac12\) is critical for \(\dot{B}^{s}_{1,\infty}\), so the smaller space is used there. Under these assumptions, the pressure, \(\nabla n^{\pm}\) and \({\bf u}^{\pm}\) satisfy \eqref{1.15} and \eqref{1.16}. The fraction densities obey the slower bound \eqref{1.17}. The exponent \(-\frac14+\frac s2\) is negative for \(s<\frac12\), zero at \(s=\frac12\), and positive for \(s>\frac12\). Thus \eqref{1.17} decays when \(s<\frac12\), is uniform at \(s=\frac12\), and for \(s>\frac12\) is only an upper bound of growth.

For the proof of Theorem \ref{3mainth}, the upper bound \eqref{1.17} does not show that the solution leaves a neighborhood of the equilibrium. Instability requires a lower bound of the same order, for data that remain small in the sense of Theorems \ref{1mainth} and \ref{2mainth}. Taking \(\frac12<s\le 1\) and initial data satisfying \eqref{1.18}, one obtains \eqref{1.19} for all sufficiently large \(t\). The right-hand side tends to infinity because \(s>\frac12\), and the exponent agrees with \eqref{1.17}. The steady state \((1,{\bf 0},1,{\bf 0})\) is therefore globally algebraically unstable, and the rate in \eqref{1.17} is sharp.

%%%%%%%%%%%%%%%%%%%%%%%%%%%%%%%%%%%%%%%%%%%%%%%
\section{\leftline {\bf{Spectral analysis and linear $L^2$ estimates.}}}
\setcounter{equation}{0}
%%%%%%%%%%%%%%%%%%%%%%%%%%%%%%%%%%%%%%%%%%%%%%%
\subsection{Reformulation} In this subsection, we first reformulate the nonlinear Cauchy problem \eqref{1.8}--\eqref{1.9} near the steady state
$(1,  {\bf 0},1,  {\bf 0})$ as follows
\begin{equation}\label{2.1}
\left\{\begin{array}{l}
\partial_{t} n^++\operatorname{div}{\bf u}^+=F_1, \\
\partial_{t}{\bf u}^{+}+\beta_1\nabla n^++\beta_2\nabla
n^--\nu^+_1\Delta {\bf u}^+-\nu^+_2\nabla\operatorname{div} {\bf u}^+-\sigma^+\nabla\Delta n^+=F_2, \\
\partial_{t} n^-+\operatorname{div}{\bf u}^-=F_3, \\
\partial_{t}{\bf u}^{-}+\beta_3\nabla n^++\beta_4\nabla
n^--\nu^-_1\Delta {\bf u}^--\nu^-_2\nabla\operatorname{div} {\bf u}^--\sigma^-\nabla\Delta n^-=F_4, \\
(n^+, {\bf u}^+, n^-, {\bf u}^-)({\bf x},0)=(n^+_0, {\bf u}^+_0, n^-_0, {\bf u}^-_0)({\bf x})\rightarrow
(0, {\bf 0}, 0, {\bf0}),~~~\hbox{as}~~~|{\bf x}|\rightarrow +\infty,\\
\end{array}\right.
\end{equation}
where $\nu_{1}^{\pm}=\frac{\mu^{\pm}}{\bar{\rho}^{\pm}}$,
$\nu_{2}^{\pm}=\frac{\mu^{\pm}+\lambda^{\pm}}{\bar{\rho}^{\pm}}>0$,
$\beta_{1}=\frac{\mathcal{C}^{2}(1,1)
\bar{\rho}^{-}}{\bar{\rho}^{+}}$,
$\beta_{2}=\beta_{3}=\mathcal{C}^{2}(1,1)$,
$\beta_{4}=\frac{\mathcal{C}^{2}(1,1)
\bar{\rho}^{+}}{\bar{\rho}^{-}}$ (which imply
$\beta_{1}\beta_{4}=\beta_{2}\beta_{3}=\beta_{2}^2=\beta_{3}^2$),
 and
the nonlinear terms are given by
\begin{align}
\label{2.2}F_{1}=&-\operatorname{div}\left(n^{+} {\bf u}^{+}\right), \\
F_{2}^{i}=&-\sqrt{\frac{\rho^-}{\rho^+}}\mathcal{C}^2\left(\left(\sqrt{\frac{\rho^-}{\rho^+}}-\beta^+\right)\partial_i n^++\left(\sqrt{\frac{\rho^+}{\rho^-}}-\beta^-\right)\partial_i n^-\right)\nonumber\\
&-\left(\sqrt{\frac{\rho^-}{\rho^+}}\mathcal{C}^2-\beta_2\beta^+\right)\partial_i \theta
-\left({\bf u}^{+} \cdot \nabla\right) u_{i}^{+} \nonumber\\
&+\mu^{+} h_{+}\left(n^{+}, n^{-}\right) \partial_{j}
n^{+} \partial_{j} u_{i}^{+}+\mu^{+} k_{+}\left(n^{+}, n^{-}\right) \partial_{j} n^{-} \partial_{j} u_{i}^{+} \nonumber\\
\label{2.3}&+\mu^{+} h_{+}\left(n^{+}, n^{-}\right) \partial_{j}
n^{+} \partial_{i} u_{j}^{+}+\mu^{+} k_{+}\left(n^{+}, n^{-}\right)
\partial_{j} n^{-}
 \partial_{i} u_{j}^{+}\\
&+\lambda^{+} h_{+}\left(n^{+}, n^{-}\right) \partial_{i} n^{+}
\partial_{j} u_{j}^{+}+\lambda^{+} k_{+}\left(n^{+},
n^{-}\right) \partial_{i} n^{-} \partial_{j} u_{j}^{+} \nonumber\\
&+\mu^{+} l_{+}\left(\rho^{+}\right) \partial_{j}^{2}
u_{i}^{+}+\left(\mu^{+}+\lambda^{+}\right) l_{+}\left(\rho^{+}\right) \partial_{i}
 \partial_{j} u_{j}^{+}, \nonumber\\
\label{2.4}F_{3}=&-\operatorname{div}\left(n^{-} {\bf u}^{-}\right), \\
F_{4}^{i}=&-\sqrt{\frac{\rho^+}{\rho^-}}\mathcal{C}^2\left(\left(\sqrt{\frac{\rho^+}{\rho^-}}-\beta^-\right)\partial_i n^++\left(\sqrt{\frac{\rho^-}{\rho^+}}-\beta^+\right)\partial_i n^-\right)\nonumber\\
&-\left(\sqrt{\frac{\rho^+}{\rho^-}}\mathcal{C}^2-\beta_3\beta^-\right)\partial_i \theta-\left({\bf u}^{-} \cdot \nabla\right) u_{i}^{-}\nonumber \\
&+\mu^{-} h_{-}\left(n^{+}, n^{-}\right) \partial_{j} n^{+}
\partial_{j} u_{i}^{-}+\mu^{-} k_{-}\left(n^{+}, n^{-}\right)
\partial_{j} n^{-} \partial_{j} u_{i}^{-} \nonumber\\
\label{2.5}&+\mu^{-} h_{-}\left(n^{+}, n^{-}\right) \partial_{j}
n^{+}
\partial_{i} u_{j}^{-}+\mu^{-} k_{-}\left(n^{+}, n^{-}\right)
 \partial_{j} n^{-} \partial_{i} u_{j}^{-} \\
&+\lambda^{-} h_{-}\left(n^{+}, n^{-}\right) \partial_{i} n^{+}
\partial_{j} u_{j}^{-}+\lambda^{-} k_{-}\left(n^{+}, n^{-}\right)
 \partial_{i} n^{-} \partial_{j} u_{j}^{-}\nonumber \\
&+\mu^{-} l_{-}\left(\rho^{-}\right) \partial_{j}^{2}
u_{i}^{-}+\left(\mu^{-}+\lambda^{-}\right) l_{-}\left(\rho^{-}\right) \partial_{i} \partial_{j} u_{j}^{-},\nonumber
\end{align}
where
\begin{equation}\label{2.8}
\left\{\begin{array}{l}
h_{+}\left(n^{+}, n^{-}\right)=\frac{\left(\mathcal{C}^{2}\alpha^{-}\right)\left(n^{+}+1, n^{-}+1\right)}{(n^++1)s_{-}^{2}\left(n^{+}+1, n^{-}+1\right)},\\
h_{-}\left(n^{+}, n^{-}\right)=-\frac{\left(\mathcal{C}^{2}
\right)\left(n^{+}+1,
n^{-}+1\right)}{(\rho^-s_{-}^{2})\left(n^{+}+1, n^{-}+1\right)},
\end{array}\right.
\end{equation}

\begin{equation}\label{2.9}
\left\{\begin{array}{l}
k_{+}\left(n^{+}, n^{-}\right)=-\frac{\mathcal{C}^{2}\left(n^{+}+1, n^{-}+1\right)}{(n^++1)(s_{+}^{2}\rho^+)
\left(n^{+}+1, n^{-}+1\right)},\\
k_{-}\left(n^{+}, n^{-}\right)=-\frac{\left(\alpha^+\mathcal{C}^{2}\right)\left(n^{+}+1, n^{-}+1\right)}{(n^-+1)
s_{+}^{2}\left(n^{+}+1, n^{-}+1\right)},\\
\end{array}\right.
\end{equation}

\begin{equation}\label{2.10}
l_{\pm}\left(\rho^\pm\right)=\frac{1}{\rho^{\pm}}-\frac{1}{\bar\rho^{\pm}}.
\end{equation}

\bigbreak

\par Define
$\tilde{U}=(\tilde{n}^+,\tilde{{\bf u}}^+,\tilde{n}^-,\tilde{{\bf u}}^-)^t$. In
terms of the semigroup theory for evolutionary equation, we will
investigate the following initial value problem for the
corresponding linear system of \eqref{2.1}:
\begin{equation}
\begin{cases}
\tilde{U}_t=\mathcal B\tilde{U},\\
\tilde{U}\big|_{t=0}={\tilde U}_0=(\tilde n^+_0, \tilde u^+_0,
\tilde n^-_0, \tilde u^-_0),
\end{cases}   \label{2.11}
\end{equation}
where the operator $\mathcal B$ is given by
\begin{equation}\nonumber\mathcal B=\begin{pmatrix}
0&-\text{div}&0&0\\
-\beta_1\nabla+\sigma^+ \nabla\Delta&\nu_{1}^{+} \Delta+\nu_{2}^{+} \nabla \otimes \nabla&-\beta_2\nabla&0\\
0&0&0&-\text{div}\\
-\beta_3\nabla&0&-\beta_4\nabla+\sigma^- \nabla\Delta&\nu_{1}^{-}
\Delta+\nu_{2}^{-} \nabla \otimes \nabla
\end{pmatrix}.\end{equation}

The solution semigroup $\text{e}^{t\mathcal B}$ has the following  decay-in-time properties, which can be found in \cite{WYZ2}.
\begin{Proposition}\label{Proposition2.1} Let $k\ge 0$, $0\le s\le 1$  and $1\le p\le 2$. Then for any $t\ge
1$, the global solution
$\tilde{{U}}=(\tilde{n}^+,\tilde{{\bf u}}^+,\tilde{n}^-,\tilde{{\bf u}}^-)^t$ of
the Cauchy problem \eqref{2.11} satisfies
\begin{equation}\label{2.12}\norm{\Lambda^{k+1}\tilde{n}^{{\pm,l}}(t)}_{L^2}\lesssim(1+t)^{-\frac{3}{2}\left(\frac{1}{p}-\frac{1}{2}\right)
-\frac{k+1}{2}}\norm{
\tilde n^\pm_0}_{L^p}+(1+t)^{-\frac{3}{2}\left(\frac{1}{p}-\frac{1}{2}\right)-\frac{k}{2}}\norm{
\tilde {\bf u}^\pm_0}_{L^p},\end{equation}
\begin{equation}\label{2.121}\norm{\Lambda^{k+1}\tilde{n}^{\pm,l}(t)}_{L^2}\lesssim(1+t)^{-\frac{3}{4}-\frac{k+1-s}{2}}\norm{
\tilde n^\pm_0}_{\dot B^s_{1,\infty}}+(1+t)^{-\frac{1}{4}-\frac{k+1-s}{2}}\norm{
\tilde {\bf u}^\pm_0}_{\dot B^s_{1,\infty}},\end{equation}
\begin{equation}\label{2.13}\norm{ \Lambda^k\left(\tilde\theta^l,\tilde{{\bf u}}^{\pm,l}\right)}_{L^2}\lesssim
(1+t)^{-\frac{3}{2}\left(\frac{1}{p}-\frac{1}{2}\right)-\frac{k}{2}}\left\|
\tilde U^\pm_0\right\|_{L^p},\end{equation}
and
\begin{equation}\label{2.131}\norm{\Lambda^k\left(\tilde\theta^l,\tilde{{\bf u}}^{\pm,l}\right)}_{L^2}\lesssim
(1+t)^{-\frac{3}{4}-\frac{k-s}{2}}\norm{
\tilde U^\pm_0}_{\dot B^s_{1,\infty}}.\end{equation}
Moreover, if $0\le s<\frac{1}{2}$,
then it holds that
\begin{equation}\label{2.14}\begin{split}\norm{\tilde{n}^{\pm,l}(t)}_{L^2} \lesssim (1+t)^{-\frac{3}{4}+\frac{s}{2}}\norm{
\tilde n^\pm_0}_{\dot B^s_{1,\infty}}+(1+t)^{-\frac{1}{4}+\frac{s}{2}}\norm{
\tilde {\bf u}^\pm_0}_{\dot B^s_{1,\infty}}.\end{split}\end{equation}
If $s=\frac{1}{2}$, it holds that
\begin{equation}\label{2.15}\begin{split}\norm{\tilde{n}^{\pm,l}(t)}_{L^2} \lesssim \norm{
\tilde U_0}_{\dot B^{-1}_{2,2}}\lesssim \norm{
\tilde U_0}_{\dot B^{\frac{1}{2}}_{1,1}}.\end{split}\end{equation}
Finally, if $0\le s\le 1$ and the initial data satisfy \eqref{1.18}, the following lower-bounds on convergence
rate hold
\begin{equation}\min\left\{\norm{\Lambda^\frac{3}{4}\tilde{n}^{+,l}(t)}_{L^2},\norm{\Lambda^\frac{3}{4}{\tilde{n}}^{-,l}(t)}_{L^2}\right\}
\geq C_0c_0(1+t)^{-\frac{5}{8}+\frac{s}{2}}.\label{2.16}
\end{equation}
\end{Proposition}

%%%%%%%%%%%%%%%%%%%%%%%%%%%%%%%%%%%%%%%%%%%%%%%
\section{\leftline {\bf{Proof of Theorem \ref{1mainth}.}}}
\setcounter{equation}{0}
%%%%%%%%%%%%%%%%%%%%%%%%%%%%%%%%%%%%%%%%%%%%%%%
By a classic argument, the global existence of solutions will be
obtained by combining the local existence result with a priori
estimates.
To begin with, we state the following local existence theorem, which can be proved by standard iteration arguments.
We refer to \cite{HaLi1,HaLi2} for details.

\begin{Proposition}[Local existence]\label{local}Let $n_{0}^{\pm}\in H^3(\mathbb{R}^3)$ and
	${\bf u}_{0}^{\pm}\in H^2(\mathbb{R}^3)$, such that
	\begin{equation*}
	\inf\limits_{x\in \mathbb{R}^3}\big{\{}n_{0}^{+}+1,n_{0}^{-}+1\big{\}}>0.
	\end{equation*}
	Then there exists a constant $T_0>0$ depending on  $\norm{n_0^\pm}_{H^3}+\norm{{\bf u}_0^\pm}_{H^2}$ such that the Cauchy problem  \eqref{2.1} admits a unique solution $\left(n^{+}, {\bf u}^{+},n^{-}, {\bf u}^{-}\right)$ on $[0,T_0)$,
	satisfying
	\begin{equation*}\inf\limits_{x\in \mathbb{R}^3, 0\leq  t\leq T_0}\big{\{}n^{+}+1,n^{-}+1\big{\}}>0,
	\end{equation*}
	and
	\[
	\begin{array}{l}
	n^{\pm} \in C^{0}\left([0, T_0] ;
	H^{3}\left(\mathbb{R}^{3}\right)\right) \cap C^{1}\left([0, T_0]
	;
	H^{2}\left(\mathbb{R}^{3}\right)\right), \\
	{\bf u}^\pm \in C^{0}\left([0, T_0] ;
	H^{2}\left(\mathbb{R}^{3}\right)\right) \cap C^{1}\left([0, T_0]
	; L^2\left(\mathbb{R}^{3}\right)\right).
	\end{array}
	\]
\end{Proposition}

 In this section, we devote ourselves to deriving the a priori energy estimates for the nonlinear system \eqref{2.1}.   We
assume  a priori that for sufficiently small $\delta>0$,
\begin{equation}\label{3.1}\norm{\left(\nabla n^\pm,P-\bar P,{\bf u}^\pm\right)}_{H^2}\le
\delta,
\end{equation}
for any $t\in [0,T_{\max})$, where $T_{\max}$ denotes the maximal time of existence the strong solution.

The norms \(\norm{n^{\pm}}_{L^2}\) are not uniformly bounded in time (see Theorem 1.3 in \cite{WYZ2}). In order to circumvent the terms containing these norms, and motivated by \cite{Bresch1}, we introduce the potential energy densities
\begin{equation}\label{3.2}G^\pm(\rho^\pm)\triangleq \rho^\pm  \int_{\bar\rho^\pm}^{\rho^\pm}\frac{P^\pm(s)-\bar P}{s^2}\mathrm{d}s,
\end{equation}
where $\bar P=P^\pm(1,1)$.  By \eqref{3.1} and Sobolev inequality, we arrive at
\begin{align}\label{3.3}\frac{1}{2}\le n^\pm+1\le 2, \frac{\bar \rho^\pm}{2}\le \rho^\pm\le 2\bar\rho^\pm.
\end{align}
Hence, we immediately conclude that there exists a positive constant $C$ such that
\begin{equation} \frac{(\rho^\pm-\bar\rho^\pm)^2}{C}\le G^\pm(\rho^\pm)\le C(\rho^\pm-\bar\rho^\pm)^2.
\end{equation}

Now, we start with  the following standard energy estimate for $(\rho^+,{\bf u}^+,\rho^-, {\bf u}^-)$.
\begin{Lemma}\label{Lemma3.2} Under the conditions of Theorem \ref{1mainth}, it holds that
\begin{align}\label{3.5}\begin{split}
\frac{\mathrm{d}}{\mathrm{d}t}&\int\left(\frac{1}{2}R^+|{\bf u}^+|^2+\frac{1}{2}R^-|{\bf u}^-|^2+\left(R^+G(\rho^+)+ R^-G(\rho^-)\right)+\frac{\sigma^+}{2}|\nabla R^+|^2+\frac{\sigma^-}{2}|\nabla R^-|^2\right)\mathrm{d}{\bf x}\\&+\int\alpha^+\left(2\mu^+|\mathbb D
{\bf u}^+|^2+\lambda^+(\operatorname{div}{\bf u}^+)^2\right)\mathrm{d}{\bf x}+\int\alpha^-\left(2\mu^-|\mathbb D
{\bf u}^-|^2+\lambda^-(\operatorname{div}{\bf u}^-)^2\right)\mathrm{d}{\bf x}=0,
\end{split}\end{align}
where $\mathbb D
{\bf u}=\frac{1}{2}\left(\nabla {\bf u}+(\nabla {\bf u})^T\right)$.
\end{Lemma}
\begin{proof}
Multiplying \eqref{1.10}$_2$ and \eqref{1.10}$_4$ by ${\bf u}^+$ and ${\bf u}^-$, respectively, integrating the resultant equation over $\mathbb R^3$, we have
\begin{align}\label{3.6}\begin{split}
\frac{1}{2}\frac{\mathrm{d}}{\mathrm{d}t}&\int\left(R^+|{\bf u}^+|^2+R^-|{\bf u}^-|^2+\sigma^+|\nabla n^+|^2+\sigma^-|\nabla n^-|^2\right)\mathrm{d}{\bf x}+\int \alpha^+\nabla P^+(\rho^+)\cdot{\bf u}^+\mathrm{d}{\bf x}\\
&+\int\alpha^+\left(\mu^+|\mathbb D
{\bf u}^+|^2+\lambda^+(\operatorname{div}{\bf u}^+)^2\right)\mathrm{d}{\bf x}+\int \alpha^-\nabla P^-(\rho^-)\cdot{\bf u}^-\mathrm{d}{\bf x}\\
&+\int\alpha^-\left(\mu^+|\mathbb D
{\bf u}^-|^2+\lambda^-(\operatorname{div}{\bf u}^-)^2\right)\mathrm{d}{\bf x}=0.
\end{split}\end{align}
Now we claim that
\begin{align}\label{3.7}\begin{split}
\int \alpha^+\nabla P^+(\rho^+)\cdot{\bf u}^+\mathrm{d}{\bf x}+\int \alpha^-\nabla P^-(\rho^-)\cdot{\bf u}^-\mathrm{d}{\bf x}
=\frac{\mathrm{d}}{\mathrm{d}t}\int \left(\alpha^+G(\rho^+)+ \alpha^-G(\rho^-)\right)\mathrm{d}{\bf x}.
\end{split}\end{align}
Assume this claim for the moment, then substituting \eqref{3.7} into \eqref{3.6} yields \eqref{3.5}. Next, we return to the proof of \eqref{3.7}. Note that
\begin{align}\nonumber\frac{\mathrm{d}}{\mathrm{d}\rho^\pm}G^\pm(\rho^\pm)=\frac{G^\pm(\rho^\pm)+P^\pm(\rho^\pm)-\bar P}{\rho^\pm}.
\end{align}
Thus we see that
\begin{align}\label{3.8}\begin{split}
\frac{\mathrm{d}}{\mathrm{d}t}&\int \alpha^+G(\rho^+)\mathrm{d}{\bf x}+\int (P^+(\rho^+)-\bar P)(\alpha^+)_t\mathrm{d}{\bf x}
 \\ &=\int \frac{G^+(\rho^+)+P^+(\rho^+)-\bar P}{\rho^+}\rho^+(\alpha^+)_t\mathrm{d}{\bf x} +\int \frac{\mathrm{d}}{\mathrm{d}\rho^+} G(\rho^+)\alpha^+(\rho^+)_t\mathrm{d}{\bf x}\\
&=\int(\alpha^+\rho^+)_t\frac{\mathrm{d}}{\mathrm{d}\rho^+} G(\rho^+)\mathrm{d}{\bf x}\\
&=\int\operatorname{div}(\alpha^+\rho^+{\bf u}^+)\frac{\mathrm{d}}{\mathrm{d}\rho^+} G(\rho^+)\mathrm{d}{\bf x}\\
&=\int \alpha^+\nabla P^+(\rho^+)\cdot{\bf u}^+\mathrm{d}{\bf x}.
\end{split}\end{align}
Similarly, we have
\begin{align}\label{3.9}\begin{split}
\frac{\mathrm{d}}{\mathrm{d}t}\int \alpha^-G(\rho^-)\mathrm{d}{\bf x}+\int (P^-(\rho^-)-\bar P)(\alpha^-)_t\mathrm{d}{\bf x}
 =\int \alpha^-\nabla P^-(\rho^-)\cdot{\bf u}^-\mathrm{d}{\bf x}.
\end{split}\end{align}
Adding \eqref{3.8} to \eqref{3.9} and using the relations $\alpha^++\alpha^-=1$ and $P^+(\rho^+)=P^-(\rho^-)$ yield \eqref{3.7}. Thus, we complete the proof of Lemma \ref{Lemma3.2}.
\end{proof}

Next, we devote ourselves to the higher-order energy estimates of \((\nabla n^{\pm},P-\bar P,{\bf u}^{\pm})\).

\smallskip
\begin{Lemma}\label{Lemma3.3} Assume that the notations and hypotheses of Theorem \eqref{1mainth}  and \eqref{3.1}  are in force,
then  it holds that
\begin{align}\label{3.10}\begin{split} &\frac{\rm d}{{\rm d}t}\left\{\beta_2\norm{\nabla \theta}_{H^{1}}^2+  \norm{\nabla  {\bf u}^+}_{H^{1}}^2+\sigma^+
 \norm{\nabla^2
 n^+}_{H^{1}}^2+\norm{\nabla {\bf u}^-}_{H^{1}}^2+\sigma^- \norm{\nabla^{2} n^-}_{H^{1}}^2
\right\}\\
&\quad +\displaystyle\nu^+_1\norm{\nabla^{2}
{\bf u}^+}_{H^{1}}^2+\nu^+_2\norm{\nabla \text{\rm div}
{\bf u}^+}_{H^{1}}^2+\nu^-_1\norm{\nabla^{2}
{\bf u}^-}_{H^{1}}^2+\nu^-_2\norm{\nabla \text{\rm div} {\bf u}^-}_{H^{1}}^2\\
&\le C\delta\left(\norm {\nabla\theta}_{H^1}^2+\norm{\nabla{\bf u}^\pm}_{L^2}^2+\norm{\nabla^3 n^\pm}_{H^1}^2\right),
\end{split}\end{align}
where $\theta=\beta^+n^++\beta^-n^-$,
$\displaystyle\beta^+=\sqrt{\frac{\beta_1}{\beta_2}}$ and
$\displaystyle\beta^-=\sqrt{\frac{\beta_4}{\beta_3}}$.
\end{Lemma}
\begin{proof}Adding $\beta^+\times\eqref{2.1}_1$ to $\beta^-\times\eqref{2.1}_3$ yields
\begin{align}\label{3.11}\partial_t\theta+\beta^+\operatorname{div}{\bf u}^++\beta^-\operatorname{div}{\bf u}^-=\beta^+F_1+\beta^-F_3.
\end{align}
For $k=1,2$,  multiplying $\nabla^k\eqref{3.11}$,  $\nabla^k\eqref{2.1}_2$ and $\nabla^k\eqref{2.1}_4$ by $\beta_2\nabla^k \theta$, $\nabla^k{\bf u}^+$ and  $\nabla^k{\bf u}^-$ respectively, summing up and then integrating the resultant equation over $\mathbb R^3$ by parts, we see that
\begin{align}\label{3.12}\begin{split}\frac{1}{2} &\frac{\rm d}{{\rm d}t}\left\{\beta_2\norm{\nabla^k \theta}_{L^{2}}^2+  \norm{\nabla^k  {\bf u}^+}_{L^{2}}^2+\sigma^+
 \norm{\nabla^{k+1}
 n^+}_{L^{2}}^2+\norm{\nabla^k {\bf u}^-}_{L^{2}}^2+\sigma^- \norm{\nabla^{k+1} n^-}_{L^{2}}^2
\right\}\\
&\quad +\displaystyle\nu^+_1\norm{\nabla^{k+1}
{\bf u}^+}_{H^{1}}^2+\nu^+_2\norm{\nabla^k \text{\rm div}
{\bf u}^+}_{H^{1}}^2+\nu^-_1\norm{\nabla^{k+1}
{\bf u}^-}_{H^{1}}^2+\nu^-_2\norm{\nabla^k \text{\rm div} {\bf u}^-}_{H^{1}}^2\\
&=\int \beta_2\nabla^k\theta\nabla^k \left(\beta^+F_1+\beta^-F_3\right)\mathrm{d}{\bf x}+\int \nabla^k F_2\nabla^k {\bf u}^+\mathrm{d}{\bf x}+\int \nabla^k F_4\nabla^k {\bf u}^-\mathrm{d}{\bf x}\\
&\quad -\sigma^+\int \nabla^{k}\Delta n^+\nabla^{k}F_1 \mathrm{d}{\bf x} -\sigma^-\int \nabla^{k}\Delta n^-\nabla^{k}F_3 \mathrm{d}{\bf x}\\
&=\sum\limits_{i=1}^5{I_i}.
\end{split}\end{align}

We shall estimate each term in the right hand side of \eqref{3.12}. By H$\ddot{\text{o}}$lder inequality, Lemmas \eqref{LemmaA.1}--\eqref{LemmaA.2}, and \eqref{3.1}, we obtain
\begin{align}\label{3.13}\begin{split}|I_1|&\le C\norm{\nabla^k\theta}_{L^2}\left(\norm{\nabla^k F_1}_{L^2}+\norm{\nabla F_2}_{L^2}\right)\\& \le C\norm {\nabla^k\theta}_{L^2}\left(\norm{n^\pm}_{L^\infty}\norm{\nabla^{k+1}{\bf u}^\pm}_{L^2}+\norm{{\bf u}^\pm}_{L^\infty}\norm{\nabla^{k+1}n^\pm}_{L^2}\right)\\& \le C\delta\left(\norm {\nabla^k\theta}_{L^2}^2+\norm{\nabla{\bf u}^\pm}_{H^2}^2\right).
\end{split}\end{align}
It follows from the expression of the pressure differential that
\begin{align}\begin{split}\label{3.14}
\nabla P&=\mathcal{C}^2\big(\rho^-\nabla n^++\rho^+\nabla n^-\big)\\
&=\mathcal{C}^2\big(\rho^--\bar{\rho}^-\big)\nabla
n^++\mathcal{C}^2\big(\rho^+-\bar{\rho}^+\big)\nabla n^-\\
&\quad+\big(\mathcal{C}^2-\mathcal{C}^2(1,1)\big)\big(\bar{\rho}^-\nabla
n^++\bar{\rho}^+\nabla
n^-\big)\\
&\quad+\mathcal{C}^2(1,1)\big(\bar{\rho}^-\nabla
n^++\bar{\rho}^+\nabla n^-\big).
\end{split}\end{align}
This together with \eqref{3.1}, \eqref{3.3} and Lemmas \eqref{LemmaA.1}--\eqref{LemmaA.2} implies  that
\begin{align}\nonumber\begin{split}\norm{\nabla P}_{L^2}\le &C \left(\norm{\mathcal{C}^2}_{L^\infty}\norm{\rho^\pm-\bar\rho^\pm}_{L^\infty}\norm{\nabla n^\pm}_{L^2}+\norm{\mathcal{C}^2-\mathcal{C}^2(1,1)}_{L^\infty}\norm{\nabla\theta}_{L^2}+\norm{\nabla\theta}_{L^2}\right)\\
\le &C\left(\delta\norm{\nabla P}_{H^1}+\norm{\nabla\theta}_{L^2}\right)
\end{split}\end{align}
and
\begin{align}\nonumber\begin{split}\norm{\nabla^2 P}_{L^2}\le& C\left(\norm{\mathcal{C}^2}_{L^\infty}\norm{\rho^\pm-\bar\rho^\pm}_{L^\infty}\norm{\nabla^2 n^\pm}_{L^2}+\norm{\nabla\left(\mathcal{C}^2\big(\rho^\pm-\bar{\rho}^\pm\big)\right)}_{L^3}\norm{\nabla n^\pm}_{L^6}\right.
\\ &\left.+\norm{\mathcal{C}^2-\mathcal{C}^2(1,1)}_{L^\infty}\norm{\nabla^2\theta}_{L^2}
+\norm{\nabla\mathcal{C}^2}_{L^3}\norm{\nabla\theta}_{L^6}+\norm{\nabla^2\theta}_{L^2}\right)\\
\le &C\left(C\delta\norm{\nabla P}_{H^1}+C\delta\norm{\nabla\theta}_{H^1}+\norm{\nabla^2\theta}_{L^2}\right).
\end{split}\end{align}
Thus, performing the similar procedure in  the above two estimates, it is easy to check that there exists a positive constant $C$ such that
\begin{align}\frac{1}{C}\norm{\nabla P}_{H^\ell}\le \norm{\nabla\theta}_{H^\ell}\le C\norm{\nabla P}_{H^\ell},\ \ell=0,1.\label{3.15}
\end{align}
Integrating by parts,  performing the similar arguments in \eqref{3.13} and making use of \eqref{3.15}, we conclude that
\begin{align}\label{3.16}\begin{split}|I_2|&=\left|\int \nabla^{k-1} F_2\nabla^{k+1}{\bf u}^+\mathrm{d}{\bf x} \right|\\
&\le C\norm{\nabla^2{\bf u}}_{H^1}\Big(\norm{\nabla P}_{H^1}\norm{\nabla n^\pm}_{H^1}+\norm{\nabla \theta}_{H^1}\norm{\nabla n^\pm}_{H^1}+\norm{\nabla{\bf u}^+}_{H^1}^2\\
&\quad+\norm{\nabla n^\pm}_{H^1}\norm{\nabla{\bf u}^+}_{H^2}+\norm{\nabla\rho^+}_{H^1}\norm{\nabla^2 {\bf u}^+}_{H^1}\Big)\\
&\le C\delta\left(\norm{\nabla\theta}_{H^1}^2+\norm{\nabla{\bf u}^+}_{H^2}^2\right).
\end{split}\end{align}
Similarly, we have
\begin{align}\label{3.17}\begin{split}|I_3|\le C\delta\left(\norm{\nabla\theta}_{H^1}^2+\norm{\nabla{\bf u}^-}_{H^2}^2\right).
\end{split}\end{align}
For the terms $I_4$ and $I_5$, it holds that
\begin{align}\label{3.18}\begin{split}|I_4|+|I_5|\le &C\norm{\nabla^{k+2}n^\pm}_{L^2}\left(\norm{n^\pm}_{L^\infty}\norm{\nabla^{k+1} {\bf u}^\pm}_{L^2}+\norm{{\bf u}^\pm}_{L^\infty}\norm{\nabla^{k+1} n^\pm}_{L^2}\right)\\
\le & C\delta\left(\norm{\nabla^{k+2}n^\pm}_{L^2}^2+\norm{\nabla {\bf u}^\pm}_{H^2}^2\right).
\end{split}\end{align}
Substituting \eqref{3.13} and \eqref{3.16}--\eqref{1.18} into \eqref{3.12} and summing $k$ from 1 to 2 yield \eqref{3.10} if $\delta$ is sufficiently small. Thus, we complete the proof of Lemma \ref{Lemma3.3}.
\end{proof}

Now, we shall establish the dissipation estimate for $n^+$ and $n^-$.
\begin{Lemma}\label{Lemma3.4} Assume that the notations and hypotheses of Theorem \eqref{1mainth}  and \eqref{3.1}  are in force,
then, for $k=0,1,2$,  it holds that
\begin{align}\label{3.19}\begin{split}\frac{\mathrm{d}}{\mathrm{d}t}&\int \left(\nabla^k{\bf u}^+\cdot\nabla\nabla^{k}n^++\nabla^k{\bf u}^-\cdot\nabla\nabla^{k}n^-\right)\mathrm{d}{\bf x}+\beta_2\norm{\nabla^{k+1}\theta}_{L^2}^2\\
&+\sigma^+\norm{\nabla^{k+2}n^+}_{L^2}^2
+\sigma^-\norm{\nabla^{k+2}n^-}_{L^2}^2\le C\left(\norm{\nabla{\bf u}^\pm}_{H^2}^2+\delta\norm{\nabla\theta}_{H^1}^2+\delta\norm{\nabla^2 n^\pm}_{H^2}^2\right).
\end{split}\end{align}
\end{Lemma}
\begin{proof} Let $k=0,1,2$. Multiplying $\nabla^k\eqref{2.1}_2$ and $\nabla^k\eqref{2.1}_4$ by $\nabla^{k}\nabla n^+$ and $\nabla^{k}\nabla n^-$ respectively, we have
\begin{align}\label{3.20}\begin{split}\beta_2&\norm{\nabla^{k+1}\theta}_{L^2}^2+\sigma^+\norm{\nabla^{k+2}n^+}_{L^2}^2
+\sigma^-\norm{\nabla^{k+2}n^-}_{L^2}^2\\
&=-\int\nabla^k\partial_t{\bf u}^+\cdot\nabla\nabla^{k}n^+\mathrm{d}{\bf x}+\int \nabla^k F_2\nabla^{k+1}n^+\mathrm{d}{\bf x}-\int\nabla^k\partial_t{\bf u}^-\cdot\nabla\nabla^{k}n^-\mathrm{d}{\bf x}+\int \nabla^k F_4\nabla^{k+1}n^-\mathrm{d}{\bf x}\\
&=\sum\limits_{i=1}^4 J_i.
\end{split}\end{align}
To handle the first term $J_1$, integrating by parts for both the $t$- and {\bf x}-variables and making use of the  continuity equation \eqref{2.1}$_1$, we obtain
\begin{align}\label{3.21}\begin{split}J_1=&-\frac{\mathrm{d}}{\mathrm{d}t}\int \nabla^k{\bf u}^+\cdot\nabla\nabla^{k}n^+\mathrm{d}{\bf x}-\int\nabla^k\operatorname{div}{\bf u}^+\nabla^kn^+_t\mathrm{d}{\bf x}\\
=&-\frac{\mathrm{d}}{\mathrm{d}t}\int \nabla^k{\bf u}^+\cdot\nabla\nabla^{k}n^+\mathrm{d}{\bf x}+\norm{\nabla^k\operatorname{div}{\bf u}^+}_{L^2}^2+\int\nabla^k\operatorname{div}{\bf u}^+\nabla^k \operatorname{div}\left(n^+{\bf u}^+\right)\mathrm{d}{\bf x}\\
\le &-\frac{\mathrm{d}}{\mathrm{d}t}\int \nabla^k{\bf u}^+\cdot\nabla\nabla^{k}n^+\mathrm{d}{\bf x}+\norm{\nabla^k\operatorname{div}{\bf u}^+}_{L^2}^2\\
&+C\left(\norm{n^+}_{L^\infty}\norm{\nabla^{k+1}{\bf u}^+}_{L^2}^2+\norm{{\bf u}^+}_{L^\infty}\norm{\nabla^{k+1}{\bf n}^+}_{L^2}\norm{\nabla^k\operatorname{div}{\bf u}^+}_{L^2}\right)\\
\le &-\frac{\mathrm{d}}{\mathrm{d}t}\int \nabla^k{\bf u}^+\cdot\nabla\nabla^{k}n^+\mathrm{d}{\bf x}+C\norm{\nabla{\bf u}^+}_{H^2}^2+C\delta\norm{\nabla^2 n^+}_{H^1}.\\
\end{split}\end{align}
Similarly, we have
\begin{align}\label{3.22}\begin{split}J_2\leq-\frac{\mathrm{d}}{\mathrm{d}t}\int \nabla^k{\bf u}^-\cdot\nabla\nabla^{k}n^-\mathrm{d}{\bf x}+C\norm{\nabla{\bf u}^-}_{H^2}^2+C\delta\norm{\nabla^2 n^-}_{H^1}.\\
\end{split}\end{align}
On the other hand, integrating by parts and recalling the estimates of $I_2$ in \eqref{3.16}, we can bound $J_2$ and $J_4$ as follows
\begin{align}\label{3.23}\begin{split}|J_2|+|J_4|&=\left|\int \nabla^{k-1} F_2\nabla^{k+2}n^+\mathrm{d}{\bf x}\right|+\left|\int \nabla^{k-1} F_4\nabla^{k+2}n^-\mathrm{d}{\bf x}\right|\\
&\le C\delta\left(\norm{\nabla\theta}_{H^1}^2+\norm{\nabla{\bf u}^\pm}_{H^2}^2+\norm{\nabla^{2} n^\pm}_{H^2}^2\right).
\end{split}\end{align}
Plugging  \eqref{3.21}--\eqref{3.23} into \eqref{3.20} and summing from $k=0$ to $2$ yield \eqref{3.19}. Thus, we complete the proof of Lemma \ref{Lemma3.4}.
\end{proof}

Now, we are in a position to prove Theorem \ref{1mainth}.

\bigskip

\noindent\textbf{Proof of Theorem \ref{1mainth}.}   Let $D>0$ be a large but fixed constant. We define the first energy functional $H_1\left(\nabla n^\pm,P-\bar P,{\bf u}^\pm\right)$ as follows
\begin{align}\begin{split}&H_1\left(\nabla n^\pm,P-\bar P,{\bf u}^\pm\right)
\\ &\quad=D\int\left(\frac{1}{2}R^+|{\bf u}^+|^2+\frac{1}{2}R^-|{\bf u}^-|^2+\left(R^+G(\rho^+)+ R^-G(\rho^-)\right)+\frac{\sigma^+}{2}|\nabla R^+|^2+\frac{\sigma^-}{2}|\nabla R^-|^2\right)\mathrm{d}{\bf x}\\
&\quad\quad +D\left\{\beta_2\norm{\nabla \theta}_{H^{1}}^2+  \norm{\nabla  {\bf u}^+}_{H^{1}}^2+\sigma^+
 \norm{\nabla^2
 n^+}_{H^{1}}^2+\norm{\nabla {\bf u}^-}_{H^{1}}^2+\sigma^- \norm{\nabla^{2} n^-}_{H^{1}}^2
\right\}\\
&\quad\quad +\sum\limits_{k=0}^2\int \left(\nabla^k{\bf u}^+\cdot\nabla\nabla^{k}n^++\nabla^k{\bf u}^-\cdot\nabla\nabla^{k}n^-\right)\mathrm{d}{\bf x}.
\end{split}\end{align}
Then, in view of \eqref{3.1} and \eqref{3.15}, it is easy to check that the energy functional $H_1\left(\nabla n^\pm,P-\bar P,{\bf u}^\pm\right)$ is equivalent to $\norm{\left(\nabla n^\pm,P-\bar P,{\bf u}^\pm\right)}_{H^2}^2$ if $D$ is suitably large. Furthermore, summing \eqref{3.19} from $k=0$ to $2$,  then adding the resultant inequality to  $D\times \left(\eqref{3.5}+\eqref{3.10}\right)$, we can conclude that there exists a positive constant $C_1$ such that
\begin{align}\label{3.25}\frac{\mathrm{d}}{\mathrm{d}t}H_1(t)+ C_1\left(\norm{\nabla\theta}_{H^2}^2+\norm{\nabla {\bf u}^\pm}_{H^2}^2+\norm{\nabla^2 n^\pm}_{H^2}^2 \right)\le 0.
\end{align}
Integrating \eqref{3.25} from $0$ to $t$ yields \eqref{1.14}. Moreover, multiplying \eqref{2.1}$_1$ by $n^+$, and then integrating the resultant equation over $\mathbb R^3$ by parts, we have
\begin{align}\nonumber\begin{split}\frac{1}{2}\frac{\mathrm{d}}{\mathrm{d}t}\norm{n^+}_{L^2}^2=&-\int n^+\operatorname{div}{\bf u}^+\mathrm{d}{\bf x}-\int n^+{\bf u}^+\cdot\nabla n^+\mathrm{d}{\bf x}-\int |n^+|^2\operatorname{div}{\bf u}^+\mathrm{d}{\bf x}\\
\le &C\left(\norm{n^+}_{L^2}\norm{\operatorname{div}{\bf u}^+}_{L^2}+\norm{n^+}_{L^\infty}\norm{n^+}_{L^2}\norm{\operatorname{div}{\bf u}^+}_{L^2}\right)\\
\le &C\norm{n^+}_{L^2}\norm{\operatorname{div}{\bf u}^+}_{L^2}.
\end{split}\end{align}
This along with \eqref{1.14} implies that
\begin{align}\label{3.26}\begin{split}\norm{n^+(t)}_{L^2}\le &\norm{n^+_0}_{L^2}+ C \int_0^t\norm{\operatorname{div}{\bf u}^+(\tau)}_{L^2}\mathrm{d}\tau\\
\le &\norm{n^+_0}_{L^2}+ C \left(\int_0^t\norm{\operatorname{div}{\bf u}^+(\tau)}_{L^2}^2\mathrm{d}\tau\right)^\frac{1}{2}\left(\int_0^t 1\mathrm{d}\tau\right)^\frac{1}{2}\\
\le & CK_0(1+t)^\frac{1}{2},
\end{split}\end{align}
and similarly, we have
\begin{align}\label{3.27}\begin{split}\norm{n^-(t)}_{L^2}
\le  CK_0(1+t)^\frac{1}{2}.
\end{split}\end{align}
Combining \eqref{3.26} and \eqref{3.27} immediately gives \eqref{1.13}. Thus, we complete the proof of Theorem \ref{1mainth}.  $\hfill\square$

%%%%%%%%%%%%%%%%%%%%%%%%%%%%%%%%%%%%%%%%%%%%%%%
\section{\leftline {\bf{Proof of Theorem \ref{2mainth}.}}}
\setcounter{equation}{0}
%%%%%%%%%%%%%%%%%%%%%%%%%%%%%%%%%%%%%%%%%%%%%%%
In this section, we devote ourselves to the long-time behavior in Theorem \ref{2mainth}. To begin with, we establish an a priori time-decay estimate for the low-frequency part \(\bigl(n^{+,l},{\bf u}^{+,l},n^{-,l},{\bf u}^{-,l}\bigr)\). The argument is based on the decay estimates for the linearized system in Proposition \ref{Proposition2.1}. More precisely, we have the following lemma.
\begin{Lemma}\label{Lemma4.1} Let  $(n^+,{\bf u}^+,n^-,{\bf u}^-)$ be a strong solution of the Cauchy problem \eqref{2.1}, then, for any $t\ge 1$, $\left( n^{+,l},{\bf u}^{+,l},n^{-,l},{\bf u}^{-,l}\right)$ satisfies the following estimates
\begin{align}\label{4.1}\begin{split}&\norm{\left(\nabla n^{\pm,l}, \theta^l,{\bf u}^{\pm,l}\right)(t)}_{L^2}\\
&\quad \le C \bigg((1+t)^{-\frac{3}{4}+\frac{s}{2}}N_0+\int_0^t(1+\tau)^{-\frac{3}{4}} \left(\norm{n^\pm}_{L^2}\norm{\nabla{\bf u}^\pm}_{L^2}+\norm{\nabla n^\pm}_{L^2}\norm{{\bf u}^\pm}_{L^2}\right.\\
&\quad \quad +\norm{\nabla n^\pm}_{L^2}\norm{P-\bar P}_{L^2}+\norm{n^\pm}_{L^2}\norm{\nabla\theta}_{L^2}+\norm{{\bf u}^\pm}_{L^2}\norm{\nabla{\bf u}^\pm}_{L^2}\\
&\quad\quad\left.+\norm{\nabla n^\pm}_{L^2}\norm{\nabla{\bf u}^\pm}_{L^2}+\norm{P-\bar P}_{L^2}\norm{\nabla^2{\bf u}^\pm}_{L^2}\right)(t-\tau)
\mathrm{d}\tau\bigg),
\end{split}\end{align}
\begin{align}\label{4.2}\begin{split}&\norm{\left(\nabla^2 n^{\pm,l}, \nabla\theta^l,\nabla{\bf u}^{\pm,l}\right)(t)}_{L^2}\\
&\quad \le C \bigg( (1+t)^{-\frac{5}{4}+\frac{s}{2}}N_0+\int_0^\frac{t}{2}(1+\tau)^{-\frac{3}{4}}\left(\norm{\nabla n^\pm}_{L^2}\norm{\left(\nabla{\bf u}^\pm,\nabla\theta\right)}_{L^2}\right)(t-\tau)
\mathrm{d}\tau\\
&\quad\quad+\int_{\frac{t}{2}}^t(1+\tau)^{-\frac{7}{4}}\left(\norm{n^\pm}_{L^2}\norm{\left({\bf u}^\pm,\theta\right)}_{L^2}\right)(t-\tau)
\mathrm{d}\tau\\
&\quad \quad+\int_0^t(1+\tau)^{-\frac{5}{4}}\left(\norm{P-\bar P}_{L^2}\norm{\nabla^2{\bf u}^\pm}_{L^2}+\norm{{\bf u}}_{L^2}\norm{\nabla{\bf u}^\pm}_{L^2}\right)(t-\tau)\mathrm{d}\tau\\
&\quad \quad+\int_0^\frac{t}{2}(1+\tau)^{-\frac{3}{4}}\left(\norm{\nabla^2 n^\pm}_{L^2}\norm{\left(\nabla{\bf u}^\pm,\theta\right)}_{L^2}+\norm{\nabla n^\pm}_{L^2}\norm{\left(\nabla^2{\bf u}^\pm,\nabla\theta\right)}_{L^2}\right)(t-\tau)
\mathrm{d}\tau\\
&\quad \quad+\int_\frac{t}{2}^t(1+\tau)^{-\frac{5}{4}}\left(\norm{\nabla n^\pm}_{L^2}\norm{\left(\nabla{\bf u}^\pm,\theta\right)}_{L^2}\right)(t-\tau)
\mathrm{d}\tau\bigg),
\end{split}\end{align}
and for $0\le s\le \frac{1}{2}$,
\begin{align}\label{4.3}\begin{split}\norm{ n^{\pm,l}}_{L^2}& \le C \bigg((1+t)^{-\frac{1}{4}+\frac{s}{2}}N_0+\int_0^t(1+\tau)^{-\frac{1}{4}} \left(\norm{n^\pm}_{L^2}\norm{\nabla{\bf u}^\pm}_{L^2}+\norm{\nabla n^\pm}_{L^2}\norm{{\bf u}^\pm}_{L^2}\right.\\
&\quad \quad +\norm{\nabla n^\pm}_{L^2}\norm{P-\bar P}_{L^2}+\norm{n^\pm}_{L^2}\norm{\nabla\theta}_{L^2}+\norm{{\bf u}^\pm}_{L^2}\norm{\nabla{\bf u}^\pm}_{L^2}\\
&\quad\quad\left.+\norm{\nabla n^\pm}_{L^2}\norm{\nabla{\bf u}^\pm}_{L^2}+\norm{P-\bar P}_{L^2}\norm{\nabla^2{\bf u}^\pm}_{L^2}\right)(t-\tau)
\mathrm{d}\tau\bigg).
\end{split}\end{align}
\end{Lemma}
\begin{proof} Setting ${U}=\left({n}^+,{{\bf u}}^+,{n}^-,{{\bf u}}^-\right)^t$ and $\mathcal F=\left(F_1,F_2,F_3,F_4\right)^t$,  it follows from Duhamel principle that
\begin{align}\nonumber U=\text{e}^{t\mathcal B}+\int_0^t \text{e}^{\tau\mathcal B} \mathcal F(t-\tau)\mathrm{d}\tau.
\end{align}
This along with Proposition \eqref{Proposition2.1}, Plancherel theorem, H\"older inequality and
Hausdorff--Young inequality leads to
\begin{align}\label{4.4}\begin{split}&\norm{\left(\nabla n^{\pm,l}, \theta^l,{\bf u}^{\pm,l}\right)(t)}_{L^2}\\
&\quad \le C\left( (1+t)^{-\frac{3}{4}+\frac{s}{2}}N_0+\int_0^t \norm{\text{e}^{\tau\mathcal B} \mathcal F(t-\tau)}_{L^2}\mathrm{d}\tau\right)\\
&\quad \le C \left((1+t)^{-\frac{3}{4}+\frac{s}{2}}N_0+\int_0^t(1+\tau)^{-\frac{3}{4}} \norm{\mathcal F(t-\tau)}_{L^1}\mathrm{d}\tau\right)\\
&\quad \le C \bigg((1+t)^{-\frac{3}{4}+\frac{s}{2}}N_0+\int_0^t(1+\tau)^{-\frac{3}{4}} \left(\norm{n^\pm}_{L^2}\norm{\nabla{\bf u}^\pm}_{L^2}+\norm{\nabla n^\pm}_{L^2}\norm{{\bf u}^\pm}_{L^2}\right.\\
&\quad \quad +\norm{\nabla n^\pm}_{L^2}\norm{P-\bar P}_{L^2}+\norm{n^\pm}_{L^2}\norm{\nabla\theta}_{L^2}+\norm{{\bf u}^\pm}_{L^2}\norm{\nabla{\bf u}^\pm}_{L^2}\\
&\quad\quad\left.+\norm{\nabla n^\pm}_{L^2}\norm{\nabla{\bf u}^\pm}_{L^2}+\norm{P-\bar P}_{L^2}\norm{\nabla^2{\bf u}^\pm}_{L^2}\right)(t-\tau)
\mathrm{d}\tau\bigg),
\end{split}\end{align}
as desired in \eqref{4.1}.

Next, let
\begin{align}\nonumber\mathcal F_1=-\left(\text{div}(n^+{\bf u}^+),\partial_i\left(\left(\sqrt{\frac{\rho^-}{\rho^+}}\mathcal{C}^2-\beta_2\beta^+\right) \theta\right),\text{div}(n^-{\bf u}^-),\partial_i\left(\left(\sqrt{\frac{\rho^+}{\rho^-}}\mathcal{C}^2-\beta_3\beta^-\right) \theta\right)\right)^t,
\end{align}
\begin{align}\nonumber\mathcal F_2&=\left(0,\left({\bf u}^{+} \cdot \nabla\right) u_{i}^{+}+\mu^{+} l_{+}\left(\rho^{+}\right) \partial_{j}^{2}
u_{i}^{+}+\left(\mu^{+}+\lambda^{+}\right) l_{+}\left(\rho^{+}\right) \partial_{i}
 \partial_{j} u_{j}^{+},\right.\nonumber\\
 &\quad \left.0,\left({\bf u}^{-} \cdot \nabla\right) u_{i}^{-}+\mu^{-} l_{-}\left(\rho^{-}\right) \partial_{j}^{2}
u_{i}^{-}+\left(\mu^{-}+\lambda^{-}\right) l_{-}\left(\rho^{-}\right) \partial_{i} \partial_{j} u_{j}^{-}\right)^t,\nonumber\end{align}
and
\begin{align}\nonumber\mathcal F_3=\mathcal F-\mathcal F_1-\mathcal F_2-\left(0,\partial_i\left(\sqrt{\frac{\rho^-}{\rho^+}}\mathcal{C}^2-\beta_2\beta^+\right) \theta,0,\partial_i\left(\sqrt{\frac{\rho^+}{\rho^-}}\mathcal{C}^2-\beta_3\beta^-\right) \theta\right)^t.
\end{align}
Then, by applying the same reasoning in \eqref{4.4}, we arrive at
\begin{align}\nonumber\begin{split}&\norm{\left(\nabla^2 n^{\pm,l}, \nabla\theta^l,\nabla{\bf u}^{\pm,l}\right)(t)}_{L^2}\\
&\quad \le C \bigg((1+t)^{-\frac{5}{4}+\frac{s}{2}}N_0+\int_0^t \norm{\text{e}^{\tau\mathcal B} \nabla\mathcal F(t-\tau)}_{L^2}\mathrm{d}\tau\bigg)\\
&\quad \le C \bigg((1+t)^{-\frac{5}{4}+\frac{s}{2}}N_0+\int_0^\frac{t}{2}(1+\tau)^{-\frac{3}{4}} \norm{\mathcal F_1(t-\tau)}_{L^\frac{3}{2}}\mathrm{d}\tau\\&\quad\quad+\int_\frac{t}{2}^t(1+\tau)^{-\frac{7}{4}} \norm{\Lambda^{-1}\mathcal F_1(t-\tau)}_{\dot B^0_{1,\infty}}\mathrm{d}\tau
+\int_0^t(1+\tau)^{-\frac{5}{4}} \norm{\mathcal F_2(t-\tau)}_{L^1}\mathrm{d}\tau\\&\quad\quad+\int_0^\frac{t}{2}(1+\tau)^{-\frac{3}{4}} \norm{\nabla\mathcal F_3(t-\tau)}_{L^1}\mathrm{d}\tau
+\int_\frac{t}{2}^t(1+\tau)^{-\frac{5}{4}} \norm{\mathcal F_3(t-\tau)}_{L^1}\mathrm{d}\tau\bigg)\\
&\quad \le C \bigg((1+t)^{-\frac{5}{4}+\frac{s}{2}}N_0+\int_0^\frac{t}{2}(1+\tau)^{-\frac{3}{4}}\left(\norm{n^\pm}_{L^6}\norm{\left(\nabla{\bf u}^\pm,\nabla\theta\right)}_{L^2}\right.\\
&\quad\quad\left.+\norm{\nabla n^\pm}_{L^2}\norm{\left({\bf u}^\pm,\theta\right)}_{L^6}\right)(t-\tau)
\mathrm{d}\tau+\int_{\frac{t}{2}}^t(1+\tau)^{-\frac{7}{4}}\left(\norm{n^\pm}_{L^2}\norm{\left({\bf u}^\pm,\theta\right)}_{L^2}\right)(t-\tau)
\mathrm{d}\tau\\
&\quad \quad+\int_0^t(1+\tau)^{-\frac{5}{4}}\left(\norm{P-\bar P}_{L^2}\norm{\nabla^2{\bf u}^\pm}_{L^2}+\norm{{\bf u}}_{L^2}\norm{\nabla{\bf u}^\pm}_{L^2}\right)(t-\tau)\mathrm{d}\tau\\
&\quad \quad+\int_0^\frac{t}{2}(1+\tau)^{-\frac{3}{4}}\left(\norm{\nabla^2 n^\pm}_{L^2}\norm{\left(\nabla{\bf u}^\pm,\theta\right)}_{L^2}+\norm{\nabla n^\pm}_{L^2}\norm{\left(\nabla^2{\bf u}^\pm,\nabla\theta\right)}_{L^2}\right)(t-\tau)
\mathrm{d}\tau\\
&\quad \quad+\int_\frac{t}{2}^t(1+\tau)^{-\frac{5}{4}}\left(\norm{\nabla n^\pm}_{L^2}\norm{\left(\nabla{\bf u}^\pm,\theta\right)}_{L^2}\right)(t-\tau)
\mathrm{d}\tau\bigg),
\end{split}\end{align}
as claimed in \eqref{4.2}. Here, in the third inequality, we have used the following element estimate:
\begin{align}\nonumber\norm{\Lambda^{-1}\nabla f}_{\dot B^0_{1,\infty}}\lesssim \norm{ f}_{\dot B^0_{1,\infty}}\lesssim \norm{ f}_{ L^1}.
\end{align}

Finally, for $0\le s\le \frac{1}{2}$,  by employing the identical approach presented  in \eqref{4.4}, we have
\begin{align}\nonumber\begin{split}\norm{ n^{\pm,l}(t)}_{L^2} &\le C\left( (1+t)^{-\frac{1}{4}+\frac{s}{2}}N_0+\int_0^t \norm{\text{e}^{\tau\mathcal B} \mathcal F(t-\tau)}_{L^2}\mathrm{d}\tau\right)\\
&\quad \le C \left((1+t)^{-\frac{1}{4}+\frac{s}{2}}N_0+\int_0^t(1+\tau)^{-\frac{1}{4}} \norm{\mathcal F(t-\tau)}_{L^1}\mathrm{d}\tau\right)\\
&\quad \le C \bigg((1+t)^{-\frac{1}{4}+\frac{s}{2}}N_0+\int_0^t(1+\tau)^{-\frac{1}{4}} \left(\norm{n^\pm}_{L^2}\norm{\nabla{\bf u}^\pm}_{L^2}+\norm{\nabla n^\pm}_{L^2}\norm{{\bf u}^\pm}_{L^2}\right.\\
&\quad \quad +\norm{\nabla n^\pm}_{L^2}\norm{P-\bar P}_{L^2}+\norm{n^\pm}_{L^2}\norm{\nabla\theta}_{L^2}+\norm{{\bf u}^\pm}_{L^2}\norm{\nabla{\bf u}^\pm}_{L^2}\\
&\quad\quad\left.+\norm{\nabla n^\pm}_{L^2}\norm{\nabla{\bf u}^\pm}_{L^2}+\norm{P-\bar P}_{L^2}\norm{\nabla^2{\bf u}^\pm}_{L^2}\right)
\mathrm{d}\tau\bigg),
\end{split}\end{align}
as desired in \eqref{4.3}. We complete the proof of Lemma \ref{Lemma4.1}.
\end{proof}

Now, we are in a position to prove Theorem \ref{2mainth}.

\smallskip

\noindent\textbf{Proof of Theorem \ref{2mainth}.} For any $t\ge 0$, we  define the time-weighted energy functional $\mathcal E(t)$ as follows
\begin{align}\begin{split}\nonumber\mathcal E(t)=&\sum\limits_{0\le \tau \le t}\Big((1+\tau)^{\frac{1}{4}-\frac{s}{2}}\norm{n^\pm(\tau)}_{L^2}+(1+\tau)^{\frac{3}{4}
-\frac{s}{2}}\norm{\left(\nabla n^\pm,\theta,{\bf u}^\pm\right)(\tau)}_{L^2}\\&+(1+\tau)^{\frac{5}{4}
-\frac{s}{2}}\norm{\left(\nabla^2 n^\pm,\nabla\theta,\nabla{\bf u}^\pm\right)(\tau)}_{H^1}\Big),
\end{split}\end{align}
where $\mathcal E(t)$ is non-decreasing.  We assume  a priori that for sufficiently small $\delta>0$,
\begin{equation}\label{4.5}\mathcal E(t)\le
\delta,
\end{equation}
for any $t\in [0,\infty)$. Then, by \eqref{3.14}, we see that
\begin{align}\nonumber\begin{split}\norm{P-\bar P}_{L^2}&\le C\norm{\Lambda^{-1}\nabla( P-\bar P)}_{L^2}\\
&\le C\Big(\norm{\mathcal{C}^2\big(\rho^--\bar{\rho}^-\big)\nabla
n^++\mathcal{C}^2\big(\rho^+-\bar{\rho}^+\big)\nabla n^-}_{L^\frac{6}{5}}\\
&\quad +\norm{\big(\mathcal{C}^2-\mathcal{C}^2(1,1)\big)\big(\bar{\rho}^-\nabla
n^++\bar{\rho}^+\nabla
n^-\big)}_{L^\frac{6}{5}}+\norm{\theta}_{L^2}\Big)
\\
&\le C\left(\norm{\rho^\pm-\bar\rho^\pm}_{L^2}\norm{\nabla n^\pm}_{L^3}+\norm{n^\pm}_{L^3}\norm{\nabla\theta}_{L^2}+\norm{\theta}_{L^2}\right)\\
&\le C\left(\delta\norm{\rho^\pm-\bar\rho^\pm}_{L^2}+\norm{n^\pm}_{L^2}^\frac{1}{2}\norm{\nabla n^\pm}_{L^2}^\frac{1}{2}\norm{\nabla\theta}_{L^2}+\norm{\theta}_{L^2}\right)\\
&\le C\left(\delta\norm{\rho^\pm-\bar\rho^\pm}_{L^2}+\left((1+t)^{-\frac{1}{4}+\frac{s}{2}}\mathcal E(t)(1+t)^{-\frac{3}{4}+\frac{s}{2}}\mathcal E(t)\right)^\frac{1}{2}\norm{\nabla\theta}_{L^2}+\norm{\theta}_{L^2}\right)\\
&\le C\left(\delta\norm{P-\bar P}_{L^2}+\delta\norm{\nabla\theta}_{L^2}+\norm{\theta}_{L^2}\right).
\end{split}\end{align}
This implies that
\begin{align}\label{4.6}\norm{P-\bar P}_{L^2}\le C\left(\delta\norm{\nabla\theta}_{L^2}+\norm{\theta}_{L^2}\right).
\end{align}
Similarly, by Lemmas \ref{LemmaA.3}--\ref{LemmaA.4} and \eqref{3.15}, it is easy to check that
\begin{align}\label{4.7}\norm{\left(P-\bar P\right)^l}_{L^2}\le C\left(\delta\norm{\nabla\theta}_{L^2}+\norm{\theta^l}_{L^2}\right).
\end{align}

Next, we split three steps to close the a prior estimate \eqref{4.5}.

{\bf Step 1:} For $t\ge 1$, substituting \eqref{4.1} and \eqref{4.7} into \eqref{3.25} that
\begin{align}\label{4.8}\begin{split}\frac{\mathrm{d}}{\mathrm{d}t}&H_1(t)+ C_2\left(H_1+\norm{\nabla\theta}_{H^2}^2+\norm{\nabla {\bf u}^\pm}_{H^2}^2+\norm{\nabla^2 n^\pm}_{H^2}^2 \right)\\
& \le C\norm{\left(\nabla n^{\pm,l}, \theta^l,{\bf u}^{\pm,l}\right)(t)}_{L^2}^2\\
& \le C \left((1+t)^{-\frac{3}{2}+s}N_0^2+\bigg(\int_0^t(1+\tau)^{-\frac{3}{4}}(1+t-\tau)^{-\frac{3}{2}+s}\mathcal E^2(t)\mathrm{d}\tau\bigg)^2\right)\\
& \le C (1+t)^{-\frac{3}{2}+s}\bigg(\delta_0^2+\mathcal E^4(t)\bigg).
\end{split}\end{align}
Applying Gronwall inequality to \eqref{4.8} gives
\begin{align}\nonumber\begin{split}H_1(t)\le &\text{e}^{-C_2t}H_1(0)+C\int_0^t \text{e}^{-C_2(t-\tau)}(1+\tau)^{-\frac{3}{2}+s}\left(\delta_0^2+\mathcal E^4(t)\right)\mathrm{d}\tau\\
\le &C\text{e}^{-C_2t}\delta_0^2+C(1+t)^{-\frac{3}{2}+s}\left(\delta_0^2+\mathcal E^4(t)\right).
\end{split}\end{align}
This implies that
\begin{align}\label{4.9} \sup\limits_{0\le \tau\le t} (1+\tau)^{-\frac{3}{2}+s}H_1(\tau)\le C\left(\delta_0^2+\mathcal E^4(t)\right).
\end{align}

{\bf Step 2:} We define the second energy functional $H_2\left(\nabla n^\pm,P-\bar P,{\bf u}^\pm\right)$ as follows
\begin{align}\label{4.511}\begin{split}&H_2\left(\nabla^2 n^\pm,\nabla P,\nabla{\bf u}^\pm\right)
\\ &\quad=D\left\{\beta_2\norm{\nabla \theta}_{H^{1}}^2+  \norm{\nabla  {\bf u}^+}_{H^{1}}^2+\sigma^+
 \norm{\nabla^2
 n^+}_{H^{1}}^2+\norm{\nabla {\bf u}^-}_{H^{1}}^2+\sigma^- \norm{\nabla^{2} n^-}_{H^{1}}^2
\right\}\\
&\quad\quad +\sum\limits_{k=1}^2\int \left(\nabla^k{\bf u}^+\cdot\nabla\nabla^{k}n^++\nabla^k{\bf u}^-\cdot\nabla\nabla^{k}n^-\right)\mathrm{d}{\bf x}.
\end{split}\end{align}
Then, in view of \eqref{3.1} and \eqref{3.15}, it is easy to check that the energy functional $H_2\left(\nabla^2 n^\pm,\nabla P,\nabla{\bf u}^\pm\right)$ is equivalent to $\norm{\left(\nabla^2 n^\pm,\nabla\theta,\nabla{\bf u}^\pm\right)}_{H^1}^2$ if $D$ is suitably large. Furthermore, summing \eqref{3.19} from $k=1$ to $2$,  then adding the resultant inequality to  $D\times \left(\eqref{3.5}+\eqref{3.10}\right)$, we can conclude that there exists a positive constant $C_3$ such that
\begin{align}\label{4.11}\begin{split}\frac{\mathrm{d}}{\mathrm{d}t}&H_2(t)+ C_3\left(H_2(t)+\norm{\nabla^2\theta}_{H^1}^2+\norm{\nabla^2{\bf u}^\pm}_{H^1}^2+\norm{\nabla^3 n^\pm}_{H^1}^2 \right)\\
&\le C\left(\norm{\nabla^2 n^{\pm,l}}_{L^2}^2+\norm{\nabla{\bf u}^{\pm,l}}_{L^2}^2+\norm{\nabla\theta^l}_{L^2}^2\right).
\end{split}\end{align}
Substituting \eqref{4.2} into \eqref{4.11}, we arrive at
\begin{align}\label{4.12}\begin{split}\frac{\mathrm{d}}{\mathrm{d}t}&H_2(t)+ C_3\left(H_2(t)+\norm{\nabla^2\theta}_{H^1}^2+\norm{\nabla^2{\bf u}^\pm}_{H^1}^2+\norm{\nabla^3 n^\pm}_{H^1}^2 \right)\\
&\le C\left(\norm{\nabla^2 n^{\pm,l}}_{L^2}^2+\norm{\nabla{\bf u}^{\pm,l}}_{L^2}^2+\norm{\nabla\theta^l}_{L^2}^2\right)\\
&\le C \Bigg((1+t)^{-\frac{5}{2}+s}N_0^2+\bigg(\int_0^\frac{t}{2}(1+\tau)^{-\frac{3}{4}}(1+t-\tau)^{-2+s}\mathcal E^2(t)\mathrm{d}\tau\bigg)^2\\
&\quad +\bigg(\int_\frac{t}{2}^t(1+\tau)^{-\frac{7}{4}}(1+t-\tau)^{-1+s}\mathcal E^2(t)\mathrm{d}\tau\bigg)^2
+\bigg(\int_0^t(1+\tau)^{-\frac{5}{4}}(1+t-\tau)^{-2+s}\mathcal E^2(t)\mathrm{d}\tau\bigg)^2\\
&\quad +\bigg(\int_0^\frac{t}{2}(1+\tau)^{-\frac{3}{4}}(1+t-\tau)^{-2+s}\mathcal E^2(t)\mathrm{d}\tau\bigg)^2
+\bigg(\int_0^t(1+\tau)^{-\frac{5}{4}}(1+t-\tau)^{-\frac{3}{2}+s}\mathcal E^2(t)\mathrm{d}\tau\bigg)^2\Bigg)\\
& \le C (1+t)^{-\frac{5}{2}+s}\bigg(\delta_0^2+\mathcal E^4(t)\bigg).
\end{split}\end{align}
Applying Gronwall inequality to \eqref{4.12} yields
\begin{align}\nonumber\begin{split}H_2(t)\le &\text{e}^{-C_3t}H_2(0)+C\int_0^t \text{e}^{-C_2(t-\tau)}(1+\tau)^{-\frac{5}{2}+s}\left(\delta_0^2+\mathcal E^4(t)\right)\mathrm{d}\tau\\
\le &C\text{e}^{-C_3t}\delta_0^2+C(1+t)^{-\frac{5}{2}+s}\left(\delta_0^2+\mathcal E^4(t)\right).
\end{split}\end{align}
This implies that
\begin{align}\label{4.13}\sup\limits_{0\le \tau\le t} (1+\tau)^{-\frac{5}{2}+s}H_2(\tau)\le C\left(\delta_0^2+\mathcal E^4(t)\right).
\end{align}

{\bf Step 3:} For $0\le s\le \frac{1}{2}$, by \eqref{4.3}, we have
\begin{align}\nonumber\begin{split}\norm{ n^{\pm,l}(t)}_{L^2}& \le C \Bigg((1+t)^{-\frac{1}{4}+\frac{s}{2}}N_0+ \int_0^t(1+\tau)^{-\frac{5}{4}}(1+t-\tau)^{-\frac{3}{2}+s}\mathcal E^2(t)\mathrm{d}\tau\Bigg)\\
&\le C(1+t)^{-\frac{1}{4}+\frac{s}{2}}\Bigg(\delta_0+\mathcal E^2(t)\Bigg).
\end{split}\end{align}
This implies that
\begin{align}\label{4.14}\begin{split}\norm{n^\pm(t)}_{L^2}&\le \norm{n^{\pm,l}(t)}_{L^2}+\norm{n^{\pm,h}(t)}_{L^2}\\
&\le C(1+t)^{-\frac{1}{4}+\frac{s}{2}}\Bigg(\delta_0+\mathcal E^2(t)\Bigg)+C(1+t)^{-\frac{1}{4}+\frac{s}{2}}
\sup\limits_{0\le\tau\le t}(1+\tau)^{\frac{3}{4}-\frac{s}{2}}\norm{\nabla n^\pm(\tau)}_{L^2}.
\end{split}\end{align}
For $\frac{1}{2}<s\le 1$, performing the similar procedure as in \eqref{3.26}, we arrive at
\begin{align}\label{4.15}\begin{split}\norm{n^\pm(t)}_{L^2}\le &\norm{n^\pm_0}_{L^2}+ C \int_0^t\norm{\operatorname{div}{\bf u}^+(\tau)}_{L^2}\mathrm{d}\tau\\
\le &\norm{n^+_0}_{L^2}+ C\int_0^t(1+\tau)^{-\frac{5}{4}+\frac{s}{2}}\mathrm{d}\tau\sup\limits_{0\le \tau\le t}(1+\tau)^{\frac{5}{4}-\frac{s}{2}}\norm{\nabla{\bf u}^\pm(\tau)}_{L^2}\\
\le & C\delta_0+C(1+t)^{-\frac{1}{4}+\frac{s}{2}}\sup\limits_{0\le \tau\le t}(1+\tau)^{\frac{5}{4}-\frac{s}{2}}\norm{\nabla{\bf u}^\pm(\tau)}_{L^2}.
\end{split}\end{align}
Combining \eqref{4.14} and \eqref{4.15} shows that
\begin{align}\label{4.16}\begin{split}\sup\limits_{0\le \tau\le t}&(1+\tau)^{\frac{1}{2}-s}\norm{n^\pm(t)}_{L^2}^2\\
&\le  C\delta_0^2+C\sup\limits_{0\le \tau\le t}(1+\tau)^{\frac{5}{2}-s}\norm{\nabla{\bf u}^\pm(\tau)}_{L^2}^2+\sup\limits_{0\le\tau\le t}(1+\tau)^{\frac{3}{2}-s}\norm{\nabla n^\pm(\tau)}_{L^2}^2.
\end{split}\end{align}
Summing $D\times\left(\eqref{4.9}+\eqref{4.13}\right)+\eqref{4.16}$ up yields
\begin{align}\nonumber\mathcal E^2(t)\le C\left(\delta_0^2+\mathcal E^4(t)\right).
\end{align}
This implies
$$\mathcal E^2(t)\le C\delta_0.$$
Therefore, we complete the proof of Theorem \ref{2mainth}.$\hfill\square$

%%%%%%%%%%%%%%%%%%%%%%%%%%%%%%%%%%%%%%%%%%%%%%%
\section{\leftline {\bf{Proof of Theorem \ref{3mainth}.}}}
\setcounter{equation}{0}
%%%%%%%%%%%%%%%%%%%%%%%%%%%%%%%%%%%%%%%%%%%%%%%

In this section, we prove the global algebraic instability of the fraction densities for \(\frac12<s\le 1\). To begin with, carrying out the same argument as for \eqref{4.4}, and using \eqref{2.16} together with Theorem \ref{2mainth}, we obtain
\begin{align}\label{5.1}\begin{split}&\min\left\{\norm{\Lambda^\frac{3}{4}n^{+,l}(t)}_{L^2},
\norm{\Lambda^\frac{3}{4}n^{-,l}(t)}_{L^2}\right\}
\\ &\quad \geq C_0c_0(1+t)^{-\frac{5}{8}+\frac{s}{2}}-C\int_0^t(1+\tau)^{-\frac{5}{8}} \norm{\mathcal F(t-\tau)}_{L^1}\mathrm{d}\tau\\
&\quad \geq C_0c_0(1+t)^{-\frac{5}{8}+\frac{s}{2}}-C\int_0^t(1+\tau)^{-\frac{5}{8}} (1+t-\tau)^{-\frac{3}{2}+s}\mathrm{d}\tau\\
&\quad \geq C_0c_0(1+t)^{-\frac{5}{8}+\frac{s}{2}}-C(1+t)^{-\frac{5}{8}+\frac{s}{2}-\frac{1-s}{2}}
\\
&\quad \geq C_4(1+t)^{-\frac{5}{8}+\frac{s}{2}},
\end{split}\end{align}
if $t$ is large enough. Therefore, it follows from the above inequality and  Lemma \ref{LemmaA.3} that
\begin{equation}\norm{\Lambda^\frac{3}{4} n^{\pm}(t)}_{L^2}\ge \norm{\Lambda^\frac{3}{4} n^{+,l}(t)}_{L^2}\ge C(1+t)^{-\frac{5}{8}+\frac{s}{2}},
\end{equation}
 which together
with the interpolation inequality
\begin{equation}\nonumber\|\Lambda^\frac{3}{4}f\|_{L^2}\leq C\|f\|_{L^2}^\frac{1}{4}\|\nabla f\|_{L^2}^\frac{3}{4},\nonumber\end{equation}
and \eqref{1.17} implies that there exists a positive constant $C_5$ such that
\begin{equation}\nonumber\min\left\{\norm{{n}^{+}(t)}_{L^2},
\norm{{{n}}^{-}(t)}_{L^2}\right\}\ge C_5(1+t)^{-\frac{1}{4}+\frac{s}{2}},
\end{equation}
as desired in \eqref{1.19} if $t$ is large enough. We complete the proof of Theorem \ref{3mainth}. $\hfill\square$

\bigskip
\appendix

%%%%%%%%%%%%%%%%%%%%%%%%%%%%%%%%%%%%%%%%%%%%%%%
\section{Analytic tools}\label{1section_appendix}
%%%%%%%%%%%%%%%%%%%%%%%%%%%%%%%%%%%%%%%%%%%%%%%
We recall the Sobolev interpolation of the Gagliardo--Nirenberg inequality, cf.   \cite[pp. 125, THEOREM]{Nirenberg}.
 \begin{Lemma}\label{LemmaA.1}
 Let $0\le i, j\le k$, then we have
\begin{equation}\nonumber
\norm{\nabla^i f}_{L^p}\lesssim \norm{  \nabla^jf}_{L^q}^{1-a}\norm{ \nabla^k f}_{L^r}^a,
\end{equation}
where $a$ satisfies
\begin{equation}\nonumber a\in\left[\frac{i}{k},1\right]\quad \text{and}\quad
\frac{i}{3}-\frac{1}{p}=\left(\frac{j}{3}-\frac{1}{q}\right)(1-a)+\left(\frac{k}{3}-\frac{1}{r}\right)a.
\end{equation}
Especially, when $p=q=r=2$, we have
\begin{equation}\nonumber
 \norm{\nabla^if}_{L^2}\lesssim \norm{\nabla^jf}_{L^2}^\frac{k-i}{k-j}\norm{\nabla^kf}_{L^2}^\frac{i-j}{k-j}.
\end{equation}

\end{Lemma}
Next, to estimate the $L^p$-norm of the spatial derivatives of the
product of two functions, we shall recall the following estimate,
 cf.  \cite{Kenig}.
\begin{Lemma}\label{LemmaA.2}
For any integer $k\ge1$, we have
 \begin{equation}\nonumber
  \norm{\nabla ^k(fg)}_{L^p} \lesssim \norm{f}_{L^{p_1}}\norm{\nabla ^kg}_{L^{p_2}} +\norm{\nabla ^kf}_{L^{p_3}}\norm{g}_{L^{p_4}},
 \end{equation}
and
 \begin{equation}\nonumber
  \norm{\nabla ^k(fg)-f\nabla^kg}_{L^p} \lesssim \norm{\nabla f}_{L^{p_1}}\norm{\nabla ^{k-1}g}_{L^{p_2}} +\norm{\nabla ^kf}_{L^{p_3}}\norm{g}_{L^{p_4}},
 \end{equation}
where $p, p_1, p_{2}, p_{3}, p_{4} \in[1, \infty]$ and
$$
\frac{1}{p}=\frac{1}{p_{1}}+\frac{1}{p_{2}}=\frac{1}{p_{3}}+\frac{1}{p_{4}}.
$$
\end{Lemma}

Finally, we recall the following two lemmas concerning the estimate for the
low-frequency part and the high-frequency part of $f$.
\begin{Lemma}\label{LemmaA.3} If $f\in L^r(\mathbb R^3)$ for any $2\leq r\leq\infty$, then we have
$$\norm{f^l}_{L^r}+\norm{f^h}_{L^r}\lesssim \norm{f}_{L^r}.$$
\end{Lemma}
\begin{proof}
For $2\leq r\leq\infty$, by virtue of Young inequality for
convolutions, for the low frequency, it holds that
\begin{equation}\label{3.361}\norm{f^l}_{L^r}\lesssim \norm{\mathfrak{F}^{-1}\phi}_{L^1}\norm{f}_{L^r}\lesssim \norm{f}_{L^r},\nonumber\end{equation}
and hence
\begin{equation}\label{3.362}\norm{f^h}_{L^r}\lesssim \norm{f}_{L^r}+\norm{f^l}_{L^r}\lesssim \norm{f}_{L^r}.\nonumber\end{equation}\end{proof}

\begin{Lemma}\label{LemmaA.4} Let $f\in H^k(\mathbb R^3)$ for any integer $k\geq 2$. Then there exists a positive constant $C_0$ such that
\begin{equation}\label{3.363}\norm{\nabla^j f^h}_{L^2}\leq C_0 \norm{\nabla^{j+1}f}_{L^2},\nonumber\end{equation}
and
\begin{equation}\nonumber\norm{\nabla^{j+1} f^l}_{L^2}\leq C_0\norm{\nabla^{j}f}_{L^2},\nonumber\end{equation}
for any $0\leq j\leq k-1$.
\end{Lemma}
\begin{proof} This lemma can be shown directly by the definitions of the low-frequency and
high-frequency of $f$ and the Plancherel theorem,  and thus we omit
the details.
\end{proof}

\section*{Acknowledgments}
Guochun Wu's research was partially supported by Natural Science Foundation of Fujian Province $\#$2026J0011524, and Natural Science Foundation of Xiamen City $\#$3502Z202673102.
Yunshun Wu' research was partially supported by National Natural Science Foundation of China $\#$12161019, $\#$12361046, Guizhou Provincial Basic Research Program (Natural Science) $\#$QKHJC-ZK [2022] YB 318, Natural Science Research Project of Guizhou Provincial Department of Education $\#$ QJJ [2023] 011. Yinghui Zhang's research was partially supported by Guangxi Natural Science Foundation $\#$2026GXNSFFA00640002 and $\#$2024GXNSFDA010071, National Natural Science
Foundation of China $\#$12271114, and the Key Laboratory of Mathematical Model and Application (Guangxi Normal University), Education Department of Guangxi Zhuang Autonomous Region.
\smallskip

{\bf A conflict of interest statement} On behalf of all authors, the corresponding author states that there is no conflict of interest.
\smallskip

{\bf Data availability statement} Data sharing is not applicable to this article as no datasets were generated or
analysed during the current study.

\end{document}